\documentclass[nohyperref]{article}

\usepackage{microtype}
\usepackage{graphicx}
\usepackage{subfigure}
\usepackage{booktabs} 

\usepackage{hyperref}

\usepackage[accepted]{icml2022}

\usepackage{amsmath}
\usepackage{amssymb}
\usepackage{mathtools}
\usepackage{amsthm}

\usepackage[capitalize,noabbrev]{cleveref}

\theoremstyle{plain}
\newtheorem{theorem}{Theorem}[section]

\newtheorem{lemma}[theorem]{Lemma}

\theoremstyle{definition}

\theoremstyle{remark}

\usepackage[textsize=tiny]{todonotes}

\icmltitlerunning{
Continuous-Time Analysis of AGM via 
Conservation Laws in Dilated Coordinate Systems
}

\usepackage{xcolor}
\usepackage{varwidth}

\newcommand{\inner}[2]{\left\langle #1 ,  #2 \right\rangle}
 
\newcommand{\norm}[1]{\left\| #1 \right\|} 

\newcommand{\br}[1]{\left[ #1 \right]}
\newcommand{\pr}[1]{\left( #1 \right)}
\newcommand{\pd}[1]{\frac{\partial}{\partial #1}}

\newcommand{\aaa}{a}
\newcommand{\bb}{b}

\newcommand{\reals}{\mathbb{R}}
\newcommand{\smu}{\sqrt{\mu}}

\begin{document}

\twocolumn[
\icmltitle{Continuous-Time Analysis of Accelerated Gradient Methods via \\
Conservation Laws in Dilated Coordinate Systems}




\begin{icmlauthorlist}
    \icmlauthor{Jaewook J. Suh}{ed}
    \icmlauthor{Gyumin Roh}{ed}
    \icmlauthor{Ernest K. Ryu}{ed}
\end{icmlauthorlist}

\icmlaffiliation{ed}{Department of Mathematical Sciences, Seoul National University, Seoul, Korea}

\icmlcorrespondingauthor{Ernest K. Ryu}{ernestryu@snu.ac.kr}

\icmlkeywords{Machine Learning, ICML}

\vskip 0.3in
]



\printAffiliationsAndNotice{}  

\begin{abstract}
We analyze continuous-time models of accelerated gradient methods through deriving conservation laws in dilated coordinate systems. Namely, instead of analyzing the dynamics of $X(t)$, we analyze the dynamics of $W(t)=t^\alpha(X(t)-X_c)$ for some $\alpha$ and $X_c$ and derive a conserved quantity, analogous to physical energy, in this dilated coordinate system. Through this methodology, we recover many known continuous-time analyses in a streamlined manner and obtain novel continuous-time analyses for OGM-G, an acceleration mechanism for efficiently reducing gradient magnitude that is distinct from that of Nesterov. Finally, we show that a semi-second-order symplectic Euler discretization in the dilated coordinate system leads to an $\mathcal{O}(1/k^2)$ rate on the standard setup of smooth convex minimization, without any further assumptions such as infinite differentiability.
\end{abstract}

\section{Introduction}
Despite the significance of acceleration within the study of first-order optimization methods, a fundamental understanding of the acceleration phenomena remains elusive. Recently, continuous-time analyses of accelerated gradient methods have been extensively pursued, even using ideas from mathematical physics. However, these continuous-time analyses still retain a component of mystery: They rely on establishing that certain energy functions are nonincreasing but do not justify the origin of such energy functions.

In this work, we present a methodology for analyzing accelerated gradient methods through deriving a conservation law, analogous to the conservation of energy of physics, in a dilated coordinate system. Namely, instead of analyzing the dynamics of $X(t)$, we analyze the dynamics of $W(t)=t^\alpha(X(t)-X_c)$ for some $\alpha\in \reals$ and $X_c\in \mathbb{R}^n$.

Through this methodology, we recover many known continuous-time analyses in a streamlined manner. Furthermore, the methodology enables us to perform a novel analysis of an ODE model of OGM-G of \citet{KimFessler2021_optimizing}, an acceleration mechanism distinct from that of \cite{Nesterov1983_method}. 
Finally, we show that a semi-second-order symplectic Euler discretization 
in the dilated coordinate system leads to an $\mathcal{O}(1/k^2)$ rate on the standard setup of smooth convex minimization, without any further assumptions such as infinite differentiability.

\subsection{Preliminaries and notation}
We review the standard definitions of convex optimization and set up the notation \cite{ Nesterov2004_introductory, BoydVandenberghe2004_convex, BauschkeCombettes2017_convex, Nesterov2018_lectures, RyuYin2022_largescale}.
Throughout the paper, we use $\reals^n $ for the underlying Euclidean space with Euclidean norm $\|\cdot\|$ and inner product $\langle\cdot ,\cdot\rangle$.
For $L>0$, $f\colon\reals^n \to \reals$ is $L$-smooth if $f$ is differentiable and
\begin{align*}
    \norm{\nabla f(x) - \nabla f(y)} \le L\norm{x-y},
    \qquad
    \forall \,x, y \in \reals^n.
\end{align*}
For $\mu>0$, $g\colon\reals^n \to \reals $ is $\mu$-strongly convex if $g(x)-(\mu/2)\norm{x}^2$ is convex.
When $f$ is differentiable and convex,
\begin{align*}
    f(x) - f(y) - \inner{\nabla f(y)}{x-y} \ge 0
\end{align*}
holds for all $x, y \in \reals^n$, and we refer to this inequality as the \emph{convexity inequality}.
Throughout this paper, consider
\begin{align}
\begin{array}{ll}
\underset{x\in \mathbb{R}^n}{\mbox{minimize}}  &f(x),
  \end{array}
  \label{eq:mainprob}
\end{align}
where $f\colon\mathbb{R}^n\rightarrow\reals$ is convex and differentiable.
When \eqref{eq:mainprob} has a minimizer, write $X_\star$ to denote a minimizer.
Write $f_\star=\inf_{x\in\mathbb{R}^n} f(x)$ for the optimal value of the problem.

\paragraph{Energy and conservation law.}
Let $A\colon (0,\infty) \to \reals$ be differentiable and $B\colon (0,\infty) \to \reals$ be integrable.
Suppose
\begin{align*}
    0 = \dot{A}(t) + B(t)
\end{align*}
holds for all $t>0$.
Then, for $0<t_0<t<\infty$, integrating from $t_0$ to $t$ 
gives us the \emph{conservation law}
\begin{align*}
    E \equiv
    A(t_0) = A(t) + \int_{t_0}^t B(s) \;ds,
\end{align*}
where the \emph{energy} $E$ is independent of time.
Moreover, if the limit $\lim_{t_0 \to 0} A(t_0)$ exists, then
\begin{align*}
    E \equiv\lim_{t_0 \to 0} A(t_0) = A(t) + \int_{0}^t B(s)\; ds.
\end{align*}

\paragraph{Partial derivatives.}
Consider a function $U(W,t)$ with variables $W=(w_1, \dots , w_n) \in \reals^n$ and $t\in\reals$. Define
\begin{align*}
    \nabla_{W}U(W,t) = \pr{ \pd{w_1} U(W,t), \dots , \pd{w_n} U(W,t) }\in \reals^n.
\end{align*}
When $W(t)$ is differentiable,
the chain rule gives us
\begin{align} \label{eq:chain rule}
    \!\!
    \frac{d}{dt} U(W(t),t) 
        &= \inner{\nabla_{W}U(W(t),t)}{\dot{W}(t)} + \pd{t} U(W(t),t).
\end{align}
To clarify, the distinction between $\frac{d}{dt}$ and $\frac{\partial}{\partial t}$ corresponds to viewing $W(t)$ as a curve dependent on $t$ or viewing $W$ as an input to $U$ independent of $t$. We clarify this notation fully in Appendix~\ref{s:pd_def}.
Then for $0<t_0<t<\infty$, integrating from $t_0$ to $t$ gives us
\begin{align}   \label{eq : Integration by parts for flow}
    &\int_{t_0}^t \inner{\nabla_{W}U(W,s)}{\dot{W}(s)}\; ds    \nonumber \\
        &= U(W(t),t) - U(W(t_0),t_0) - \int_{t_0}^t \pd{s}U(W,s) \;ds.    \nonumber
\end{align}

\subsection{Prior work}
In convex optimization and machine learning, the classical goal is to reduce the function value efficiently.
In the smooth convex setup, Nesterov's celebrated accelerated gradient method (AGM) \cite{Nesterov1983_method} achieves an accelerated rate of $\mathcal{O}(1/k^2)$.
Recently, the optimized gradient method (OGM) \cite{KimFessler2016_optimized} improved the rate of AGM by a factor of $2$, and this rate is in fact exactly optimal \cite{Drori2017_exact}.
In the smooth strongly convex setup, the strongly convex AGM (SC-AGM) \citep[2.2.22]{Nesterov2018_lectures} achieves an accelerated rate.
The review by \citet{dAspremontScieurTaylor2021_acceleration} provides a comprehensive historical review.

The study of first-order convex optimization algorithms efficiently reducing the squared gradient norm was initiated by \citet{Nesterov2012_how}.
For smooth non-convex minimization, gradient descent (GD) achieves an $\mathcal{O}((f(x_0)-f_\star)/k)$ rate \citep[Proposition~3.3.1]{Nemirovsky1999_optimization}.
In the smooth convex setup, OGM-G \cite{KimFessler2021_optimizing} achieves an $\mathcal{O}(f(x_0)-f_\star)/k^2)$ rate. 
M-OGM-G \cite{ZhouTianSoCheng2022_practical} and OBL-G${}_\flat$ \cite{ParkRyu2021_optimal} are variants of OGM-G achieving similar rates.
Combining AGM with OGM-G \citep[Remark~2.1]{Nesterov2018_lectures} yields an $\mathcal{O}(\norm{x_0-x_\star}^2/k^4)$ rate, which matches the $\Omega( \norm{x_0-x_\star}^2/k^4 )$ lower bound of \citep{Nemirovsky1991_optimality,Nemirovsky1992_informationbased} and is therefore optimal.

An ODE model for the heavy ball method with constant friction, i.e., constant damping, 
was introduced by \citet{Polyak1964_methods} and
follow-up work studying variations flourished
\cite{AttouchAlvarez2000_heavy, AlvarezAttouch2001_inertiala,
AttouchCzarnecki2002_asymptotic,
AlvarezAttouchBolteRedont2002_secondorder, AttouchBolteRedont2002_optimizing, AttouchMaingeRedont2012_secondorder,AttouchCzarnecki2017_asymptotic, BotCsetnek2017_secondorder, BotCsetnek2019_secondorder, AdlyAttouch2020_finitea, AdlyAttouchVo2021_asymptotic, AujolDossalRondepierre2021_convergence, AujolDossalRondepierre2022_convergence}.
The study of ODE models of AGM and accelerated mirror descent with vanishing damping was initiated by  \citet{SuBoydCandes2014_differential,SuBoydCandes2016_differential,KricheneBayenBartlett2015_accelerated}.
Specifically, \citet{SuBoydCandes2014_differential} studied the dynamics of $0 = \ddot{X} + \frac{r}{t} \dot{X} + \nabla f(X)$ and proved $f(X(t))-f_\star \le (r-1)^2\norm{X_0-X_\star}^2/(2t^2)$ for $r\ge3$.
\citet{AttouchChbaniPeypouquetRedont2018_fast} improved the constant of this bound for $r>3$.
For $r<3$, \citet{AttouchChbaniRiahi2019_rate} established an $\mathcal{O}(t^{-2r/3})$  rate.
Improved rates under the additional, so-called, $\mathbf{H}_1(\gamma)$ hypothesis were established by \citet{AujolDossalRondepierre2019_optimal, SebbouhDossalRondepierre2019_nesterovs, ApidopoulosAujolDossalRondepierre2021_convergence}.
A wide range of variations of the ODE with vanishing damping were also studied
\cite{AttouchChbani2015_fast,May2017_asymptotic,AttouchCabotChbaniRiahi2018_rate,AttouchChbaniRiahi2018_combining, AttouchCabot2018_convergence,AttouchChbaniRiahi2019_fasta,AttouchPeypouquet2019_convergence, AttouchLaszlo2020_newtonlike,AttouchChbaniFadiliRiahi2020_firstorder,AttouchBalhagChbaniRiahi2021_fast,AttouchFadiliKungurtsev2021_effect, AttouchLaszlo2021_continuous,AttouchCabot2017_asymptotic, AttouchLaszlo2021_convex, BotCsetnekLaszlo2021_tikhonov, AttouchBalhagChbaniRiahi2022_damped,AttouchChbaniFadiliRiahi2021_convergence}.
Similar analyses were extended to differential inclusions for non-differentiable functions
\cite{AttouchMainge2011_asymptotic, AttouchPeypouquet2016_rate, AujolDossal2017_optimala, ApidopoulosAujolDossal2017_second, ApidopoulosAujolDossal2018_differential},
monotone inclusions \cite{BotCsetnek2016_second, BotCsetnek2018_convergence, BotCsetnekLaszlo2018_secondorder, BotHulett2022_second},
primal-dual methods \cite{BotNguyen2021_improved},
and splitting methods \citet{FrancaRobinsonVidal2018_admm, Hassan-MoghaddamJovanovic2021_proximal, FrancaRobinsonVidal2021_gradient, AttouchChbaniFadiliRiahi2021_fast}.

This intense study of ODEs modeling optimization algorithms motivated the development of tools utilizing the following ideas:
variational principle and Lagrangian mechanics \cite{WibisonoWilsonJordan2016_variational, Jordan2018_dynamical, ZhangOrvietoDaneshmand2021_rethinking, WilsonRechtJordan2021_lyapunov};
duality gap and convex-analytical techniques \cite{DiakonikolasOrecchia2019_approximate};
Hamiltonian mechanics \cite{DiakonikolasJordan2021_generalized};
control theory \cite{HuLessard2017_dissipativity};
continuous-time complexity lower bounds  \cite{MuehlebachJordan2020_continuoustime}; 
and
perturbation analysis of physics, leading to the high-resolution ODE \cite{ShiDuJordanSu2021_understanding}.

The study of continuous-time models, in turn, motivated the study of discretizing such ODEs to obtain implementable algorithms.
Discretizing ODEs with vanishing damping
\cite{WibisonoWilsonJordan2016_variational,AttouchCabotChbaniRiahi2018_inertial, AttouchCabot2018_convergencea,AttouchChbaniRiahi2019_fast, AttouchChbaniFadiliRiahi2020_firstorder,AdlyAttouch2020_finite,AttouchCabot2020_convergence, AttouchChbaniRiahi2020_convergence,AdlyAttouch2021_firstorder,AdlyAttouchLe2021_first, AdlyAttouchVo2021_newtontype,DiakonikolasJordan2021_generalized} 
and discretizing alternate ODEs \cite{ScieurRouletBachdAspremont2017_integration,WilsonMackeyWibisono2019_accelerating, MuehlebachJordan2019_dynamical, ZhangSraJadbabaie2019_acceleration} 
have been studied.
Specifically, \citet{ZhangMokhtariSraJadbabaie2018_direct} achieved an $\mathcal{O}(1/k^2)$ rate using the Runge--Kutta discretization on the ODE by \citet{SuBoydCandes2014_differential} under additional assumptions.

The study of using symplectic integrators, a discretization scheme designed to conserve energy \cite{HairerLubichGerhard2006_geometric}, for discretizing the ODE models was initiated by \citet{BetancourtJordanWilson2018_symplectic} and was further developed in a series of work \cite{MaddisonPaulinTehODonoghueDoucet2018_hamiltonian, FrancaSulamRobinsonVidal2020_conformal, FrancaSulamRobinsonVidal2020_conformala,  MuehlebachJordan2021_optimization,FrancaJordanVidal2021_dissipative}.
However, these approaches did not obtain an asymptotic $\mathcal{O}(1/k^2)$ rate in the sense usually considered in optimization.
An $\mathcal{O}(1/k^2)$ rate was obtained by \citet{ShiDuSuJordan2019_acceleration} combining symplectic integration with the high-resolution ODE framework.

Recently, \citet{EvenBerthierBachFlammarionHendrikxGaillardMassoulieTaylor2021_continuized} introduced the ``continuized'' framework of accelerated gradient methods, which uses a stochastic jump process to perform randomized discretizations. The framework can utilize the simpler continuous-time analysis while producing an implementable (but randomized) discrete algorithm with rate $\mathcal{O}(1/k^2)$.

\subsection{Contribution}
The central thesis, the main contribution, of this paper is that continuous-time analyses of accelerated gradient methods significantly simplify under an alternate dilated coordinate system. We establish this claim by presenting a methodology analyzing the ODEs by deriving conservation laws in dilated coordinate systems and recovering many prior analyses in a streamlined manner. We then use the methodology to perform the first continuous-time analysis of OGM-G, whose acceleration mechanism was understood far less than the acceleration mechanism of Nesterov.

Furthermore, we show that the coordinate change can also benefit the analysis of discretizations. 
Specifically, we apply a semi-second-order symplectic Euler discretization 
in the dilated coordinate system to obtain an $\mathcal{O}(1/k^2)$ rate in the standard setup of smooth convex minimization, without any further assumptions such as infinite differentiability. 
This is the first result of its kind, in the precise sense clarified in 
Section~\ref{ss:section5discussion}, and it will be interesting to see, in future work, to what extent discretizations exploiting our dilated coordinates can achieve competitive rates.

\section{Conservation laws from dilated coordinates} \label{s:section2}
Our main methodology for continuous-time analysis is to perform a coordinate change
and then obtain a conservation law.
In this section, we quickly exhibit this methodology applied to the classical AGM ODE and then present a generalized form which we will use in later sections.

Consider problem \eqref{eq:mainprob}. Assume a minimizer of $f$ exists and write $X_\star$ for a minimizer of $f$. (We do not assume the minimizer is unique.)
Write $f_\star=f(X_\star)$.
The AGM ODE presented by \citet{SuBoydCandes2014_differential} is
\begin{equation}    \label{eq:AGMODE}
    0 = \ddot{X} + \frac{3}{t} \dot{X} + \nabla f(X)
\end{equation}
with initial condition $X(0)=X_0$, $\dot{X}(0)=0$. 
Here, $X\colon[0,\infty) \to \reals^n$ is a function of the time $t$, but we often write $X$ in place of $X(t)$ for the sake of notational brevity.
Consider the dilated coordinate $W=t^{\alpha}(X-X_\star)$ with a yet undetermined $\alpha \in \reals$.
The ODE in the $W$ coordinate is
\begin{equation}    \label{eq:ExpandingCoordinateODE}
    0 = \frac{1}{t^\alpha}\ddot{W} + \frac{3-2\alpha}{t^{\alpha+1}} \dot{W} + \nabla_{W} U(W,t)
\end{equation}
with 
\begin{equation}    \label{ftn : Potential Energy for r3}
    U(W,t) = \frac{\alpha(\alpha-2)}{2t^{\alpha+2}} \norm{W}^2
                +  t^{\alpha} \pr{ f \pr{ X(W,t) } - f_\star }
\end{equation}
and $X(W,t)=\frac{W}{t^\alpha}+X_\star$.
Since $U$ contains $t^\alpha(f(X)-f_\star)$, we choose $\alpha=2$ in anticipation of the $\mathcal{O}(1/t^2)$ rate to get
\begin{align}       \label{eq:ExpandingCoordinateODE2}
    0   &= \frac{1}{t^2}\ddot{W} - \frac{1}{t^{3}} \dot{W} + \nabla_{W} U(W,t). 
\end{align}
Taking the inner product between $\dot{W}$ and \eqref{eq:ExpandingCoordinateODE2} 
and using \eqref{eq:chain rule}, we get
\begin{align*}
    0   &= \frac{d}{dt} \pr{ \frac{1}{2t^2}\norm{\dot{W}}^2 } + \inner{\nabla_{W} U(W,t)}{\dot{W}(t)}  \\
        &= \frac{d}{dt} \pr{ \frac{1}{2t^2}\norm{\dot{W}}^2 + U(W(t),t) } - \pd{t} U(W(t),t).
\end{align*}
The corresponding conservation law is
\begin{align*}   
    &E\equiv 2\norm{X_0-X_\star}^2  \\
        &= \lim_{t_0\to0} \pr{ \frac{1}{2t_0^2}\norm{\dot{W}(t_0)}^2 + U(W(t_0),t_0) }    \\
        &= \frac{1}{2t^2} \norm{\dot{W}(t)}^2  + U(W(t),t) - \int_{0}^{t} \pd{s} U(W(s),s) \;ds.
\end{align*}
From $\pd{t}X(W,t) = -\frac{2}{t^3} W = -\frac{2}{t}(X-X_\star)$, we get
\begin{align*} 
    -\pd{t} U(W,t) 
        &= -\pd{t} t^{2} \pr{ f \pr{ X(W,t) } - f_\star }     \\
        &= 2 t \big( f_\star - f \pr{ X } - \inner{ \nabla f(X) }{ X_\star - X } \big)
\end{align*}
and
\begin{align}   \label{ftn:Energy for AGM ODE}
     E  &\equiv  2\norm{X_0-X_\star}^2    \nonumber \\ 
        &=  t^2 \pr{ f(X) - f_\star } + \frac{1}{2} \norm{t\dot{X} + 2(X-X_\star)}^2    \\&\quad \nonumber
             + \int_{0}^{t} 2s \big( f_\star - f \pr{ X } - \inner{ \nabla f(X) }{ X_\star - X } \big) \;ds  
\end{align}
for all $t\ge 0$.
Since $f$ is convex, the integrand is nonnegative, and we conclude
\begin{align*}
    f(X) - f_\star 
        \le \frac{E}{t^2} = \frac{2\norm{X_0-X_\star}^2}{t^2}.
\end{align*}

\paragraph{General form of conservation laws.}
We now generalize the previous analysis for later sections.
Let $U\colon \reals^{n}\times\reals \to \reals$, and consider the ODE
\begin{equation*}    \label{eq:CanonicalFormOfODE}
    0 = \aaa(t) \ddot{W} + \bb(t)\dot{W} + \nabla_W U(W,t).
\end{equation*}
Take the inner product with $\dot{W}$ and integrate to obtain the conservation law
\begin{align}   \label{eq : Conservation Law for ODE}
    E   &\equiv \frac{a(t_0)}{2} \norm{\dot{W}(t_0)}^2 + U(W(t_0),t_0)     \\ \nonumber
        &= \frac{\aaa(t)}{2} \norm{\dot{W}(t)}^2
        + \int_{t_0}^{t} \pr{\bb(s)-\frac{\dot{\aaa}(s)}{2}} \norm{\dot{W}(s)}^2 ds       \\ \nonumber
        &\quad + U(W(t),t)     
        - \int_{t_0}^{t} \pd{s} U(W(s),s) \;ds.
\end{align}

Note that if $\aaa(t)=1$ and $U(W,t)=U(W)$, then this convservation law is nothing but the familiar conservation of energy in physics;
within $E$, the first term $(1/2)\| \dot{W} \|^2$ is kinetic energy, the second term $\int_{t_0}^{t} \bb \| \dot{W} \|^2 ds $ is energy dissipated way as heat due to friction, the third term $U(W)$ is potential energy, and the fourth term vanishes as the potential $U$ is independent of time.

Throughout this paper, we consider dilated coordinates of the form $W=e^{\gamma(t)}(X-X_c)$ for some $X_c\in \mathbb{R}^n$.
As a consequence, $U(W,t)$ will contain $e^{\gamma(t)}( f(X(W,t)) - f(X_c) )$. The convexity inequality enters the integral of $\pd{s} U(W,s)$ through the identity
\begin{align*}   
    &-\pd{t} e^{\gamma(t)} \pr{ f \pr{ X(W,t) } - f(X_c) }     \\
        &= \dot{\gamma}(t) e^{\gamma(t)} \big( f(X_c) - f \pr{ X } - \inner{ \nabla f(X) }{ X_c - X } \big).
\end{align*}
Note, if $e^{\gamma(t)} = 1$ for all $t$, i.e. if there is no coordinate change, then $\dot{\gamma}(t)=0$ and the convexity inequality does not enter the conservation law. In this sense, the coordinate change is essential for our analysis to utilize convexity.

\paragraph{Connection with Lyapunov analyses.}
Our analyses based on conservation laws are not fundamentally different from the Lyapunov analyses of the prior work.
The first two terms of the conservation law for the AGM ODE 
    $$ \Phi(t) = t^2 \pr{ f(X) - f_\star } + \frac{1}{2} \norm{t\dot{X} + 2(X-X_\star)}^2, $$ 
form the exact Lyapunov function of \citet{SuBoydCandes2014_differential}. Once $\Phi(t)$ is stated, it is relatively straightforward to verify $\dot{\Phi}(t)\le0$ through direct differentiation.
The conservation laws of Section~\ref{s:section3} also contain Lyapunov functions of prior work \cite{AttouchChbaniRiahi2019_rate,AujolDossal2017_optimal,AujolDossalRondepierre2019_optimal}. 

The analyses of prior work often start by stating a Lyapunov function of unclear origin and then proceed with the analysis. In truth, these Lyapunov functions are obtained through many hours of trial and error. A core motivation of our work is to provide a systematic methodology for obtaining such Lyapunov functions.

The closely related prior work of \citet{DiakonikolasJordan2021_generalized} presents a methodology based on Hamiltonian mechanics.
While they also provide a unified methodology for analyzing continuous-time models of accelerated gradient methods, there are some key differences that we further clarify in Appendix~\ref{appendix : Difference between DiakonikolasJordan2021_generalized and us}.
One key difference is that while we start from a given ODE and derive conservations laws, \citet{DiakonikolasJordan2021_generalized} start from a Hamiltonian with ``potential energy`` and a ``kinetic energy'' terms and derive the ODE.
From our framework, a $\norm{W}^2$ term arises naturally as in \eqref{ftn : Potential Energy for r3} and as in the third term of \eqref{eq:EnergyGeneralizedForAlphaR}, but $\norm{W}^2$ does not arise from the approach of \citet{DiakonikolasJordan2021_generalized}.
Our analyses of the generalized AGM, SC-AGM, and OGM-G ODEs crucially rely on using the $\norm{W}^2$ term and therefore cannot be obtained by the methodology of \citet{DiakonikolasJordan2021_generalized} as is.

\section{Continuous-time analyses of Nesterov-type acceleration via conservation laws in dilated coordinate systems}
\label{s:section3}

Again, consider problem \eqref{eq:mainprob}. Assume a minimizer of $f$ exists and write $X_\star$ for a minimizer of $f$.
Write $f_\star=f(X_\star)$.
\citet{SuBoydCandes2016_differential} presented the generalized ODE 
\begin{align}   \label{eq:AGMODEWithr}
    0 = \ddot{X} + \frac{r}{t}\dot{X} + \nabla f(X)
\end{align}
and provided Lyapunov analyses for $r\ge3$.
We consider the dilated coordinate $W=t^\alpha(X-X_\star)$ and follow a similar line of reasoning as that of Section~\ref{s:section2} to obtain the conservation law
\begin{align}   \label{eq:EnergyGeneralizedForAlphaR}
    &E    \equiv t^{\alpha} \pr{ f \pr{ X } - f_\star}        \nonumber
        + \frac{1}{2} t^{\alpha-2}\norm{ t\dot{X} + \alpha (X-X_\star) }^2   \\ &\qquad 
        + \frac{\alpha(\alpha +1 - r)}{2} t^{\alpha-2} \norm{X-X_\star}^2  \\ \nonumber &\qquad 
        + \int_{t_0}^{t} \bigg( \frac{(2r-3\alpha)s^{\alpha-3}}{2} \norm{s \dot{X} + \alpha(X-X_\star)}^2 \\ \nonumber &\qquad\qquad\quad
            + \frac{\alpha(\alpha +1 - r)(\alpha+2)}{2} s^{\alpha-3} \norm{X-X_\star}^2 \bigg) \;ds \\ \nonumber &\qquad
        + \int_{t_0}^{t} \alpha s^{\alpha-1} \pr{  f_\star - f (X)
            - \inner{ \nabla f (X) } { X_\star - X } }  \;ds .
\end{align}
Note that when $r=3$, $\alpha=2$, and $t_0=0$, half of the terms vanish and the conservation law reduces to \eqref{ftn:Energy for AGM ODE}.

Throughout this section, we present the analysis results based on conservation laws while deferring the detailed derivations to Appendix~\ref{s:sectionB}.

\subsection{AGM ODE $r>3$}
\label{ss:agmrge3}
Let $r>3$. Plug $\alpha=2$ and $t_0=0$ into \eqref{eq:EnergyGeneralizedForAlphaR} and evaluate integrals as described in Appendix~\ref{appendix : r>3 detail} to get
\begin{align*}
    E 
    &\equiv (5-r) \norm{ X_0 - X_\star }^2 \\ 
    &= -2(r-3) \norm{X_0-X_\star}^2 \\ &\quad
        + t^{2} \pr{ f(X) - f_\star }
        + \frac{1}{2} \norm{ t \dot{X} + 2(X-X_\star)}^2    \\ &\quad
        + (r-3) \norm{X-X_\star}^2   
        + \int_{0}^{t} \frac{r-3}{s} \norm{s \dot{X} }^2  ds \\&\quad
        + \int_{0}^{t} 2 s \pr{  f_\star - f (X)
            - \inner{ \nabla f (X) } { X_\star - X } }  ds.
\end{align*}
All terms depending on $t$ are nonnegative when $r>3$.
Thus $E+ 2(r-3) \norm{X_0-X_\star}^2\ge t^2(f(X)-f_\star)$ holds, and we conclude
    $$ f(X)-f_\star \le \frac{(r-1)\norm{X_0-X_\star}^2}{t^2} .$$

This rate improves upon the rate  $f(X)-f_\star \le \frac{(r-1)^2\norm{X_0-X_\star}^2}{2t^2}$ by \citet{SuBoydCandes2014_differential}
and matches the rate of \citet{AttouchChbaniPeypouquetRedont2018_fast}.
This conservation law also implies $E\ge(r-3)\norm{X-X_\star}^2$, and boundedness of $\norm{X-X_\star}$ can be used to establish convergence of $X(t)$
\cite{ChambolleDossal2015_convergence, AttouchChbaniPeypouquetRedont2018_fast}.

\subsection{AGM ODE $r<3$}  \label{subsection : r<3}
\label{ss:agmrle3}
Let $0\le r<3$.
Plug $\alpha=\frac{2r}{3}$ to (\ref{eq:EnergyGeneralizedForAlphaR}) to get
\begin{align*}   
    E &= t^{\frac{2r}{3}} \pr{ f \pr{ X } - f_\star}  
        + \frac{r(3-r)}{9} t^{\frac{2r}{3}-2} \norm{X-X_\star}^2   \\ &\,\, 
        + \frac{1}{2} t^{\frac{2r}{3}-2}\norm{ t\dot{X} + \frac{2r}{3} (X-X_\star) }^2  \\  &\,\,
        + \int_{t_0}^{t} \frac{2}{27} r(3-r)(3+r) s^{\frac{2 r}{3}-3} \norm{ X - X_\star }^2 \, ds    \\ &\,\,
        + \int_{t_0}^{t} \frac{2r}{3} s^{\frac{2r}{3}-1} \pr{  f_\star - f (X)
        - \inner{ \nabla f (X) } { X_\star - X } }  \;ds .
\end{align*}
We let the starting time be nonzero, i.e., $t_0>0$,  to ensure all of the terms do not blow up.
All terms are nonnegative.
Thus $E \ge t^{\frac{2r}{3}}(f(X)-f_\star)$, and we conclude
    $$ f(X)-f_\star \le \frac{E}{t^{\frac{2r}{3}}}.$$
This recovers the result of \citet{AttouchChbaniRiahi2019_rate}.

\subsection{AGM ODE with growth condition}
\label{sec:growth condition}
\citet{AujolDossalRondepierre2019_optimal} consider convex functions satisfying the so-called ``$\mathbf{H}_1(\gamma)$ hypothesis'', defined as
\begin{align*}
    f(x)-f_\star \le \frac{1}{\gamma} \inner{\nabla f(x)}{x-X_\star},\qquad\forall x\in \mathbb{R}^n
\end{align*}
for a $\gamma\ge1$, and obtain improved rates.
To utilize the $\mathbf{H}_1(\gamma)$ hypothesis, rather than the convexity inequality, we rescale the ODE by multiplying $t^\beta$ and then obtain the conservation law \eqref{eq : Conservation Law for ODE} with the rescaled ODE.
The derivations are detailed in Appendix~\ref{appendix : AGM ODE growth condition}.
With values $\alpha = \frac{2r}{\gamma+2}$ and $\beta = \frac{2(\gamma - 1)r}{\gamma + 2}$ we get
\begin{align*}
    &E 
    \equiv t^{\frac{2\gamma r}{\gamma+2}}(f(X) - f_\star)
        + \frac{1}{2} t^{\frac{2\gamma r}{\gamma+2} - 2} \norm{ t\dot{X} + \alpha (X- X_\star) }^2       \\  &\,\,
        + \frac{r( 2 - \gamma (r-1) )}{(\gamma + 2)^2} t^{\frac{2\gamma r}{\gamma+2} - 2} \| X - X_\star\|^2 \\  &\,\,
     + \int_{t_0}^t \frac{ 2 r ( 2r + 2 - \gamma (r - 1) ) (2 - \gamma (r - 1))}{(\gamma + 2)^3}    
     \\ &\qquad\qquad\qquad\qquad\qquad\qquad\qquad
     s^{\frac{2\gamma r}{\gamma+2} - 3} \| X - X_\star \|^2 ds \\ &\,\,
    +  \int_{t_0}^t s^{\frac{2\gamma r}{\gamma+2} - 1} \frac{2\gamma r}{\gamma+2}
        \\ &\qquad\qquad\qquad
            \left(  f_\star - f(X) - \frac{1}{\gamma} \inner{ \nabla f(X) }{ X_\star - X } \right) ds.
\end{align*}

When $\gamma\ge1$ and $r \le 1+\frac{2}{\gamma}$, all terms are nonnegative, and we get 
\begin{align*}
    f(X) - f_\star \le \frac{E}{t^{\frac{2\gamma r}{\gamma+2}}},
\end{align*}
which recovers the result of \cite{AujolDossalRondepierre2019_optimal}.
Note that this rate is better than that of Section~\ref{subsection : r<3} since $\frac{2\gamma r}{\gamma+2} \ge \frac{2r}{3}$ for $\gamma \ge 1$.

\subsection{SC-AGM}
\label{ss:sc-agm}
\citet{WilsonRechtJordan2021_lyapunov} presented the following ODE of the strongly convex accelerated gradient method (SC-AGM)
\begin{align}   \label{eq:SC_AGM_ODE}
    0 = \ddot{X} + 2\sqrt{\mu} \dot{X} + \nabla f(X)
\end{align}
with initial condition $X(0)=X_0$, $\dot{X}(0)=0$, where $\mu>0$ is the strong convexity parameter of $f$.

Consider the dilated coordinate $W = e^{\sqrt{\mu} t}(X-X_\star)$. 
The resulting conservation law with $t_0=0$ is 
\begin{align*}
    E 
    &\equiv f(X_0) - f_\star  \\
    &= - \frac{\mu}{2} \norm{X_0-X_\star}^2     \\&\quad
        + e^{\smu t} \pr{ f(X)-f_\star + \frac{1}{2} \norm{ \dot{X} + \smu(X-X_\star) }^2 } \\&\quad
        + \int_{0}^{t} \frac{\smu e^{\smu s}}{2} \norm{ \dot{X}  }^2 ds    
        + \int_{0}^{t} \smu e^{\smu s} \big( ... \big) \, ds,
\end{align*}
where
\begin{align*}
    (...)
    &= f_\star - f(X) - \inner{\nabla f(X)}{X_\star-X} - \frac{\mu}{2} \norm{X-X_\star}^2  \\
    &\ge 0.
\end{align*}
The inequality follows from $\mu$-strong convexity of $f$.
All the terms depending on $t$ are nonnegative,
thus $E + \frac{\mu}{2} \norm{X_0-X_\star}^2 \ge e^{\smu t}(f(X)-f_\star)$, and we conclude
\begin{align*}
    f(X)-f_\star \le e^{-\smu t} \pr{ f(X_0)-f_\star + \frac{\mu}{2} \norm{X_0-X_\star}^2}. 
\end{align*}
This recovers the result of \cite{WilsonRechtJordan2021_lyapunov}.

\subsection{Gradient flow} \label{s:gradient flow}
We conclude this section by showing that dilated coordinates also simplify the analysis of the gradient flow ODE
\begin{align*}
    0 = \dot{X} + \nabla f(X)
\end{align*}
with $X(0)=X_0$, which is a first-order ODE model of gradient descent.

Consider the dilated coordinate $W=t(X-X_\star)$.
With $a(t)=0$ in \eqref{eq : Conservation Law for ODE},
we get the conservation law with $t_0=0$ 
\begin{align*}
    &E   \equiv -\frac{1}{2}\norm{X_0-X_\star}^2 \\
        &= t \pr{ f (X) - f_\star } + \frac{1}{2} \norm{X-X_\star}^2 - \norm{X_0-X_\star}^2 \\
        &\!\!
            +  \int_{0}^t \!\! s \norm{ \dot{X} } ^2 \!ds 
            +  \!\!\int_{0}^t\!\! \pr{ f_\star - f(X) - \inner{\nabla f(X)}{X_\star-X} }\;ds.
\end{align*}
We recover the well-known result
\begin{align*}
    f(X) - f_\star \le \frac{\norm{X_0-X_\star}^2}{2t}.
\end{align*}

\section{Continuous-time analysis of OGM-G}
\label{s:section4}
We now present a novel ODE model of OGM-G \cite{KimFessler2021_optimizing},
which optimally reduces the squared gradient magnitude (rather than the function value) for smooth convex minimization.
Consider problem \eqref{eq:mainprob}. Assume $f_\star=\inf_{x\in\mathbb{R}^n}f(x)>-\infty$.
(We do not assume a solution exists.) 
Following steps similar to those of \citet{SuBoydCandes2014_differential} with OGM-G, 
we obtain the OGM-G ODE
\begin{align*}  
    0 = \ddot{X} - \frac{3}{t-T} \dot{X} + 2\nabla f(X)
\end{align*}
for $t\in(0,T)$ with initial value $X(0)=X_0$, $\dot{X}(0)=0$.
The precise derivation of the OGM-G ODE and the calculations throughout this section are presented in Appendix~\ref{s:sectionC}.

Choose the dilated coordinate $W = (T-t)^{\alpha}(X-X_c)$ for some $X_c\in \mathbb{R}^n$.
Since we expect the rate $\mathcal{O}\pr{1/T^2}$, we choose $\alpha=-2$.
The corresponding conservation law is
\begin{align*}
    &E  \equiv \frac{2}{T^2} (f(X_0)-f(X_c)) \\ &
        =  \frac{2}{(T-t)^2}  \pr{ f(X) - f(X_c) } 
        -  \frac{2}{(T-t)^4}  \norm{X-X_c}^2  \\&  
        +  \frac{1}{2(T-t)^4}  \norm{ (T-t) \dot{X} + 2(X-X_c)}^2 \\& 
        + \int_{0}^{t} \! 
              \frac{4}{(T-s)^3} \pr{  f(X_c) - f (X) - \inner{ \nabla f (X) } { X_c-X } }  ds.    
\end{align*}

\subsection{OGM-G ODE $r=-3$}
We now establish an $\mathcal{O}(1/T^2)$ rate on $\norm{\nabla f(X(T))}^2$ via a conservation law. At first, this may seem curious as the conservation law contains no terms directly involving $\nabla f(X)$.

We first characterize the dynamics of the solution to the OGM-G ODE near the terminal time $t=T$.
\begin{lemma} \label{lemma:continuous_extension_at_T}
Let $X\colon[0,T)\rightarrow\mathbb{R}^n$ be the solution to the OGM-G ODE.
We can continuously extend $X(t)$, $\dot{X}(t)$, $\ddot{X}(t)$ to $t=T$ with
\vspace{-0.1in}
    \begin{align*}
        \dot{X}(T) &= 0,\quad
        \ddot{X}(T) = \lim_{t\to T^{-}} \frac{\dot{X}(t)}{t-T} = \nabla f(X(T)).
    \end{align*}
\end{lemma}
\begin{proof}[Proof outline]
For simplicity, assume $\lim_{t\rightarrow T^{-}}\dot{X}(t)$ and 
$\lim_{t\rightarrow T^{-}}\ddot{X}(t)=\lim_{t\rightarrow T^{-}} \frac{\dot{X}(t)-\dot{X}(T)}{t-T}$ exist.
We will formally prove these assumptions in Appendix~\ref{appendix : Property of OGM-G ODE}.

Consider the conservation law with $\alpha=0$ and $X_c=X_0$:\vspace{-0.1in}
\begin{align*}
    E     
    \equiv \frac{1}{2}\norm{\dot{X}}^2 + 2( f(X) - f(X_0) )   
    + \int_{0}^{t} \frac{3\norm{\dot{X}}^2}{T-s}\; ds.
\end{align*}
Since $E$ is independent of time and since the first two terms are bounded, we have $\int_{0}^{T} \frac{3\norm{\dot{X}}^2}{T-s} ds<\infty$. 
The finite integral implies $\lim_{t\rightarrow T^{-}}\dot{X}(t)=0$.
Furthermore,\vspace{-0.05in}
\begin{align*}
    0   &= \lim_{t\to T^{-}} \pr{ \ddot{X}(t) - \frac{3}{t-T} \dot{X}(t) + 2\nabla f(X(t)) }    \\
        &= - 2 \ddot{X}(T)  + 2\nabla f(X(T)).
\\[-0.25in]\tag*{\qedhere}
\end{align*}
\end{proof}


We now prove the promised result.
\begin{theorem} \label{thm:main_result_for_OGMG}
Let $X\colon[0,T]\rightarrow\mathbb{R}^n$ be the extended solution to the OGM-G ODE.
Then $X$ exhibits the rate
$$ \norm{\nabla f(X(T)) }^2 \le \frac{4\pr{ f(X_0)-f(X(T)) }}{T^2} 
            \le \frac{4\pr{ f(X_0)-f_\star }}{T^2}. $$
\end{theorem}
\begin{proof}
Consider the conservation law with $X_c=X(T)$ and define the Lyapunov function 
\begin{align*}
    \Phi(t) &= 
        \frac{2}{(T-t)^2}  \pr{ f(X) - f(X(T)) } \\& 
        -  \frac{2}{(T-t)^4}  \norm{X-X(T)}^2  \\&  
        +  \frac{1}{2(T-t)^4}  \norm{ (T-t) \dot{X} + 2(X-X(T))}^2.
\end{align*}
Then $\Phi(t)$ is monotonically nonincreasing by the conservation law,
and so $\Phi(0)\ge \lim_{t\to T^{-}} \Phi(t)$.

By applying L'H\^{o}pital's rule,
\begin{align*}
    \lim_{t\to T^{-}} \frac{f(X(t))-f(X(T))}{(T-t)^2}&
        = \frac{1}{2} \norm{\nabla f(X(T))}^2  \\
    \lim_{t\to T^{-}} \frac{X(t)-X(T)}{(T-t)^2}&
        = \frac{1}{2} \nabla f(X(T)).  
\end{align*}
Therefore, 
\vspace{-0.1in}
\begin{align*}
    \lim_{t\to T^{-}} \Phi(t) 
        &= \norm{\nabla f(X(T))}^2 - \frac{1}{2}\norm{\nabla f(X(T))}^2 + 0 \\
        &= \frac{1}{2}\norm{\nabla f(X(T))}^2
\end{align*}
and we conclude
\vspace{-0.05in}
\begin{gather*}
    \frac{1}{2} \norm{\nabla f(X(T)) }^2 \le \frac{2}{T^2} \pr{ f(X_0) - f(X(T)) }. \\[-0.25in]\tag*{\qedhere}
\end{gather*}
\end{proof}

In the proof of Theorem~\ref{thm:main_result_for_OGMG}, $\nabla f$ does not explicitly appear in the conservation law and only arises at the terminal time $T$ due to Lemma~\ref{lemma:continuous_extension_at_T}.
For this reason, we can establish a bound on $\norm{\nabla f(X(t))}^2$ only at the terminal time.

\citet{LeeParkRyu2021_geometric} presented the first Lyapunov analysis of the discrete-time OGM-G.
We show in Appendix~\ref{appendix : OGM-G discrete correpondence} that the Lyapunov function of Theorem~\ref{thm:main_result_for_OGMG} is the continuous-time analog of the Lyapunov function of \citet{LeeParkRyu2021_geometric}.
The discrete-time analysis for OGM-G also establish a rate on $\norm{\nabla f(x_k)}^2 $ only for the terminal iteration $k=K$.

\subsection{OGM-G ODE for $r<-3$}    \label{section : OGM-G r<-3}
Following \citet{SuBoydCandes2014_differential}, we generalize the OGM-G ODE to general $r$:\vspace{-0.1in}
\begin{align*} \label{eq:generalized_OGMG_ODE}
    0 = \ddot{X} + \frac{r}{t-T} \dot{X} + 2\nabla f(X).
\end{align*}
In Appendix~\ref{appendix : Property of OGM-G ODE}, 
we directly extend the arguments of Lemma~\ref{lemma:continuous_extension_at_T} to conclude $\lim_{t\to T^{-}} \frac{\dot{X}(t)}{t-T} = -\frac{2}{r+1}\nabla f(X(T))$.

With the dilated coordinate $W=(T-t)^{-2}(X-X(T))$, we get the conservation law
\begin{align*}
    &E  \equiv \frac{2}{T^2} ( f(X) - f(X(T)) ) + \frac{r+3}{T^4} \norm{X-X(T)}^2 \\&\,\,
        \,\,\,\!\!\!= \frac{2}{(T-t)^{2}} \pr{ f(X) - f(X(T)) }  
                + \frac{r+1}{(T-t)^{4}} \norm{X-X(T)}^2 \\&\,\,
        + \frac{1}{2(T-t)^{4}} \norm{ (T-t) \dot{X} + 2(X-X(T))}^2 \\&\,\,  
        + \int_{0}^{t} \frac{(-(r+3))}{(T-s)^{5}} \norm{(T-s) \dot{X} + 2(X-X(T))}^2 ds  \\&\,\,
        \!\!\!\!+\!\! \int_{0}^{t}\!\! \frac{4}{ (T-s)^{3} } \pr{  f(X(T)) \!-\! f (X)
            \!- \inner{ \nabla f (X) } { X(T) - X } }  ds.
\end{align*}

\begin{theorem} \label{thm:generalized_result_for_OGMG}
Let $X\colon[0,T]\rightarrow\mathbb{R}^n$ be the extended solution to the OGM-G ODE with $r<-3$.
Then,
        $$ \norm{\nabla f(X(T)) }^2 
            \le \frac{2(-1-r)\pr{ f(X_0)-f(X(T)) }}{T^2}  $$
\end{theorem}
\begin{proof}[Proof outline]
The arguments are similar to those of Theorem~\ref{thm:main_result_for_OGMG}:
Define a Lyapunov function $\Phi(t)$ based on the conservation law 
and consider the inequality $\Phi(0)\ge\lim_{t\to T^{-}} \Phi(t)$.
Details are presented in Appendix~\ref{appendix : Detail for generalized_result_for_OGMG}.
\end{proof}

\subsection{Obtaining $\norm{\nabla f(X(T))}^2 \le \mathcal{O}(1/T^4)$ with OGM $+$ OGM-G ODE}
We state a simple technique to obtain an $\mathcal{O}( \norm{x_0-x_\star}^2/T^4)$ rate from the $\mathcal{O}((f(X_0)-f_\star)/T^2)$ rate of the OGM-G ODE. This technique is based on the idea of \citet{Nesterov2012_how}, \citet{NesterovGasnikovGuminovDvurechensky2020_primaldual} to concatenate AGM with OGM-G to obtain a $\norm{\nabla f(x_K)}^2 \le \mathcal{O} ( \norm{x_0-x_\star}^2/K^4 )$ rate.

If one starts the AGM ODE with $X^\mathrm{F}(0)=X^\mathrm{F}_0$ and $\dot{X}^\mathrm{F}(0)=0$, the terminal solution $X^\mathrm{F}(T)$ satisfies $f(X^\mathrm{F}(T)) - f_\star  \le 2\norm{X_0-X_\star}^2/T^2$.
Then we start the OGM-G ODE with $X^\mathrm{G}(0)=X^\mathrm{F}(T)$ and $\dot{X}^\mathrm{G}(0)=0$ and obtain the solution $X^\mathrm{G}(T)$ satisfying $\norm{\nabla f(X^\mathrm{G}(T))}^2 \le 4 (f(X^\mathrm{G}(0)) - f_\star )/T^2$.
Concatenating these two guarantees, we obtain $\norm{\nabla f(X^\mathrm{G}(T))}^2 \le 8 \norm{X_0-X_\star}^2/T^4$.

\section{Discretization in dilated coordinates via semi-second-order symplectic Euler}
\label{s:section5}

In this section, we show that discretizing the AGM ODE ($r=3$) using a semi-second-order symplectic Euler discretization in the dilated coordinate system leads to an algorithm with an $\mathcal{O}(1/k^2)$ rate.
Despite the extensive prior work on continuous-time analyses and discretizations of the AGM ODE, obtaining an accelerated rate through a direct and ``natural'' discretization has been surprisingly tricky. 
Our result is the first to accomplish this, in the precise sense clarified in Section~\ref{ss:section5discussion}.

Again, the ODE \eqref{eq:AGMODE}, restated, is $0 = \ddot{X} + \frac{3}{t} \dot{X} + \nabla f(X)$.
With $W=t^{2}(X-X_\star)$, the ODE \eqref{eq:ExpandingCoordinateODE2}, restated, is
\begin{align}       \tag{\ref{eq:ExpandingCoordinateODE2}}
    0   &= \frac{1}{t^2}\ddot{W} - \frac{1}{t^{3}} \dot{W} + \nabla_{W} U(W,t).
\end{align}
We first identify a generalized coordinate $W$ and conjugate momentum $P$ to replace $X$ and $\dot{X}$.
The dilated coordinate $W=t^{2}(X-X_\star)$ has been chosen, so we determine the generalized momentum via the Lagrangian formulation.

Recall from \eqref{ftn : Potential Energy for r3} that
    $ U(W,t) = t^2\pr{ f(X(W,t)) - f_\star } $.
Define the Lagrangian as 
\begin{align*}
  L(W,\dot{W},t) = \frac{1}{2t} \norm{\dot{W}}^2 - t \; U(W,t).
\end{align*}
Then the Euler--Lagrange equation $\frac{d}{dt}\nabla_{\dot{W}}L = \nabla_{W}L$ yields the ODE \eqref{eq:ExpandingCoordinateODE2} and $P=\nabla_{\dot{W}}L = \frac{\dot{W}}{t} =t\dot{X}+2(X-X_\star)$ is the conjugate momentum.
Express \eqref{eq:ExpandingCoordinateODE2} in $W$ and $P$:
\begin{align*}
    \dot{P}&=-t\nabla f(X(W,t))\\
    \dot{W}&=tP     
\end{align*}
and
$\ddot{W} =P - t^2 \, \nabla f \pr{ X(W,t) }$.

Inspired by the symplectic Euler \citep{HairerLubichGerhard2006_geometric} and velocity Verlet integrators
\citep{Verlet1968_computer,SwopeAndersenBerensWilson1982_computer,AllenTildesley2017_computer}
we consider alternating updates of $W$ and $P$ but use a second-order update for $W$:
\begin{align*}   
    P(t+h)  &\approx P(t) - t\nabla f(X) h     \nonumber \\
    W(t+h)  &\approx W(t) + \dot{W}(t) h + \ddot{W}(t) \frac{h^2}{2}      \\
            = W(t) &+ t P(t) h + \big(  P(t) - t^2\nabla f(X(W,t)) \big) \frac{h^2}{2} . \nonumber
\end{align*}
We refer to this method as a semi-second-order symplectic Euler. This discretization is also an instance of the Nystr\"{o}m method \citep{HairerLubichGerhard2006_geometric}.

Identifying  $w_k$ and $p_k$ with $W(hk)$ and $P(hk)$ and defining $x_k$ through $w_k=h^2k^2(x_k-X_\star)$, we get the method
\begin{align*}
    \!\!\!p_{k+1} &= p_k - k h^2 \nabla f \pr{ x_k }     \\
    \!\!\!x_{k+1} &= \!\frac{k^2}{(k+1)^2} \pr{ \!x_{k} \!-\! \frac{h^2}{2} \nabla{f} \pr{ x_k } \!} +\!\frac{2k+1}{(k+1)^2}  \pr{ \frac{p_{k+1}}{2} \!+\!X_\star\!}\!.
\end{align*}
Finally, letting $s=h^2$, $\theta_k=\frac{k}{2}$ 
 and  $z_k = \frac{p_k}{2} + X_\star$, we get
    \begin{align}
    \label{eq:AGM ODE discretization}
        x_{k}^+ &= x_k - \frac{s}{2} \nabla f(x_k)  \nonumber \\
        z_{k+1} &= z_{k} - s\theta_k \nabla f(x_k)    \\ 
        x_{k+1} &= \frac{\theta_{k}^2}{\theta_{k+1}^2} x_{k}^{+} 
                    + \pr{1-\frac{\theta_{k}^2}{\theta_{k+1}^2}} z_{k+1}       \nonumber 
    \end{align}
for $k=0,1,\dots$. 
The starting point is $x_0=z_0=X_0\in\mathbb{R}^n$,
since $z_0$ corresponds to $\frac{P(0)}{2}+X_\star = X_0$.

\begin{theorem}     \label{thm: convergence rate for r=3}
    Assume $f$ is convex and $L$-smooth. Assume $f$ has a minimizer $X_\star$. 
    For $s\in \left( 0,\frac{2}{L} \right]$,  \eqref{eq:AGM ODE discretization} exhibits the rate \vspace{-0.05in}
        \[
        f(x_k^{+}) - f_\star \le \frac{2\norm{ X_0 - X_\star }^2}{s k^2}   .
        \vspace{-0.05in}
        \]
\end{theorem}
\begin{proof}[Proof outline]
    The proof is based on the Lyapunov analysis $\Phi_k\le \Phi_{k-1}\le\dots\le \Phi_0$ with\vspace{-0.05in}
    \begin{align*}
    \!
       \Phi_k &\!=\! 2 c_k \theta_k^2 \pr{ f(x_k) \!-\! f_\star \!-\! \frac{s}{4} \norm{\nabla f(x_k) }^2 \!} +  \frac{1}{s} \norm{ z_{k+1} - X_\star }^2 \vspace{-0.05in}
    \end{align*}
    and $c_k=\frac{\theta_{k+1}}{\theta_{k+1}^2-\theta_k^2}$
    for $k=0,1,\dots$.
    The details are presented in Appendix~\ref{appendix:proof_for_thm_5.1}.
\end{proof}

\subsection{Discussion}
\label{ss:section5discussion}

\paragraph{Hamiltonian mechanics.}
Some may wonder what can be said from a Hamiltonian mechanics perspective. We discuss this matter briefly in Appendix~\ref{s:hamiltonian}, and \cite{DiakonikolasJordan2021_generalized, FrancaJordanVidal2021_dissipative} pursues this direction deeply. Here, we point out the quick observation that the explicit time-dependence of the Lagrangian makes the Hamiltonian time-dependent, and this time-dependence makes the Hamiltonian a non-conserved quantity. Therefore, the classical theory of symplectic integrators is not immediately applicable, 
but we nevertheless use our method and obtain an accelerated rate.

\paragraph{Prior discretizations.}
The discretization of \cite{WibisonoWilsonJordan2016_variational} achieves an $\mathcal{O}(1/k^2)$ rate, but, arguably, this discretization 
``does not flow natural from the dynamical-systems framework'' \citep[p.\ 529]{Jordan2018_dynamical}.
\citet{ZhangMokhtariSraJadbabaie2018_direct} achieved an accelerated rate with a Runge--Kutta method, but their $\mathcal{O}(1/k^2)$ rate requires the additional assumption of infinite differentiability.
\citet{ShiDuSuJordan2019_acceleration} used a symplectic integrator 
with $\dot{X}$ as the momentum (no coordinate change) and achieved an $\mathcal{O}(1/k^2)$ rate, but they crucially rely on the high-resolution ODE formulation.
\citet{FrancaJordanVidal2021_dissipative} proposed a generalized symplectic integrator and established  $\mathcal{O}(1/k^2)$ rate for exponentially large $k$ depending on the stepsize, but their rate does not hold for all $k\in \mathbb{N}$.
\citet{EvenBerthierBachFlammarionHendrikxGaillardMassoulieTaylor2021_continuized} introduced alternative ``continuized'' framework and obtained $\mathcal{O}(1/k^2)$ with randomized discretizations.
On the other hand, our result is a direct, non-randomized discretization of the AGM ODE that achieves an $\mathcal{O}(1/k^2)$ rate without making additional assumptions or using a high-resolution formulation.

\paragraph{Discretized rate surpasses AGM.}
The rate of Theorem~\ref{thm: convergence rate for r=3} with $s=\frac{2}{L}$ is
\vspace{-0.1in}
\[
    f(x_k^{+}) - f_\star
        \le  \frac{L\norm{ X_0 - X_\star }^2}{k^2}.
\]
Interestingly, this rate is smaller (better) than the rate of Nesterov's AGM by a factor of 2 \cite{Nesterov1983_method} but is slightly larger (worse) than the exact optimal rate of OGM 
\citep{DroriTeboulle2014_performance,KimFessler2016_optimized,Drori2017_exact}.
This improvement seems to be in part due to the choice of Lyapunov function, inspired by \cite{ParkParkRyu2021_factorsqrt2}, that allows a tighter analysis. By taking the continuous-time limit of AGM and then discretizing, we arrived at a discretized algorithm that is \emph{better} than the original AGM.

\paragraph{Interpreting $z_k$ as conjugate momentum.}
\citet{LeeParkRyu2021_geometric} point out that many known accelerated gradient methods have an auxiliary $z_k$-sequence satisfying a geometric structure. In our analysis of the AGM ODE, we identify that $z_k$ is (up to a factor-$2$ scaling and translation with $X_\star$) the conjugate momentum $P=\dot{W}/t=t\dot{X} + 2(X - X_\star)$ of the dilated coordinate $W=t^2(X-X_\star)$. 

Moreover, we've observed that this interpretation of the $z$-variables as conjugate momenta of the dilated coordinate systems (with some rescaling and translation) also holds in other setups, including the SC-AGM and the OGM-G setups. Specifically, when we discretize the ODEs in the dilated coordinate systems $W(t)$, the discretized methods closely resemble the known accelerated methods, and the $z$-variables roughly correspond to conjugate momenta $P(t)$. We leave the formalization and development of this observation as future work.

\section{Conclusion} \label{s:section6}
This work presents a methodology for analyzing continuous-time models of accelerated gradient methods through deriving conservation laws in dilated coordinate systems. Using this methodology, we recover many known continuous-time analyses in a streamlined manner and obtain novel continuous-time analyses of OGM-G.

We hypothesize that our dilated coordinates can simplify analyses of other setups beyond those explored in Sections~\ref{s:section3} and \ref{s:section4}. For example, exploring the use of dilated coordinates in stochastic differential equations modeling stochastic optimization and investigating whether dilated coordinates generally simplify discretization, as was the case for the AGM ODE ($r=3$) in Section~\ref{s:section5}, are interesting directions of future work. Finally, finding a more fundamental understanding of the interpretation of $z_k$ as the conjugate momentum would also be interesting.

\section*{Acknowledgements}
JJS and EKR were supported by the Samsung Science and Technology Foundation grant (Project Number SSTF-BA2101-02).
We thank Jongmin Lee for valuable discussions about OGM-G. 
We thank Chanwoo Park for reviewing the manuscript and providing valuable
feedback. 
Finally, we thank the anonymous reviewers for their thoughtful comments.

\bibliography{conservation_law_agm}
\bibliographystyle{icml2022}

\newpage
\appendix
\onecolumn

\section{Partial derivative notation}
\label{s:pd_def}
For $U\colon\reals^n \times \reals \to \reals$, we assign symbols $W\in \reals^n$ and $t\in\reals$ for the inputs, i.e., we write $U(W,t)$.
At the same time, we consider the curve $W\colon \reals\rightarrow\reals^n$ a function of $t\in \reals$, i.e., we write $W(t)$.
When we provide the curve $W(t)$ as the first input to $U$, we get $U(W(t),t)$, which is now a function solely of $t\in \reals$, and we can take the total derivative $\frac{d}{dt}$ of it.
Using the chain rule of vector calculus, we get
\begin{align*}
\frac{d}{dt} U(W(t),t)&=\inner{(D_1 U)(W(t),t)}{\dot{W}(t)}+ (D_2U)(W(t),t)
\end{align*}
where $D_1 U$ is the derivative of $U(\cdot,\cdot) $ with respect to the first $n$ coordinates
and $D_2 U$ is the derivative of $U(\cdot,\cdot) $ with respect to the last coordinate.
When $U(W,t)$ is viewed as a function of $W$ and $t$ (when $W$ is an input variable independent of $t$ rather than a curve), then
\[
D_1 U=\nabla_W U,\qquad D_2 U=\pd{t} U.
\]
Write $\nabla_{W}U(W(t),t)$ to mean take the partial derivative of $U(W,t)$ with respect to $W$ and then plug in $W=W(t)$.
Likewise, write $\pd{t} U(W(t),t)$ to mean take the partial derivative of $U(W,t)$ with respect to $t$ and then plug in $W=W(t)$.
Finally, we can write
\begin{align*}
\frac{d}{dt} U(W(t),t)&=\inner{(D_1 U)(W(t),t)}{\dot{W}(t)}+ (D_2U)(W(t),t)\\
&= \inner{\nabla_{W}U(W(t),t)}{\dot{W}(t)} + \pd{t} U(W(t),t).
\end{align*}

\section{Comparison with \cite{DiakonikolasJordan2021_generalized}} \label{appendix : Difference between DiakonikolasJordan2021_generalized and us}

\citet{DiakonikolasJordan2021_generalized} present a methodology based on Hamiltonian mechanics, and their goal is also to provide a unified methodology for analyzing continuous-time models of accelerated gradient methods.
However, our methodology differs from that of \citet{DiakonikolasJordan2021_generalized} in the following three ways.
\begin{itemize}
    \item We start from a given ODE and derive conservations laws, while \citet{DiakonikolasJordan2021_generalized} start from a Hamiltonian and derive the ODE.
    \item In our framework, different choices of `$\alpha$' produce different conservation laws for one fixed ODE, but in \cite{DiakonikolasJordan2021_generalized} different choices of `$\alpha$' corresponds to different ODEs and different corresponding energies.
    \item Our framework accommodates translation with respect to an arbitrary ``center point'' $X_c$.
\end{itemize}
Our analyses of the AGM, SC-AGM, and OGM-G ODEs crucially rely on these differences and therefore cannot be obtained by the methodology of \citet{DiakonikolasJordan2021_generalized} as-is:
\begin{itemize}
    \item The approach of \citet{DiakonikolasJordan2021_generalized} does not lead to a Lyapunov function or a conservation law containing $\norm{W}^2$. Many of our results crucially rely on using an energy $U(W,t)$ with the $\norm{W}^2$ term.
    \item The translation with respect to $X_c=X(T)$ is essential for the analysis of OGM-G ODE in Theorem~\ref{thm:main_result_for_OGMG}.
\end{itemize}

\section{Omitted calculations of Section~\ref{s:section3}}
\label{s:sectionB}

\subsection{Conservation law for generalized $r$} \label{appendix : derivation of mechanical energy for general r}
We start with ODE \eqref{eq:AGMODEWithr}
\begin{align*}
    0 = \ddot{X} + \frac{r}{t}\dot{X} + \nabla f(X) .
\end{align*}
Now consider the coordinate change $W=t^{\alpha}(X-X_\star)$.
Then we see
\begin{align*}
    W       &= t^\alpha(X-X_\star)     \\
    \dot{W} &= t^\alpha\dot{X} + \alpha t^{\alpha-1}(X-X_\star)  \\ 
    \ddot{W} &= t^\alpha\ddot{X} + 2\alpha t^{\alpha-1}\dot{X} + \alpha(\alpha-1)t^{\alpha-2}(X-X_\star).  
\end{align*}
From this, we can rewrite $X$, $\dot{X}$, $\ddot{X}$ in terms of $W$, $\dot{W}$, $\ddot{W}$,
\begin{align*}
    X &= \frac{W}{t^\alpha} + X_\star     \\
    \dot{X} &= \frac{\dot{W}}{t^\alpha} - \alpha\frac{W}{t^{\alpha+1}}     \\
    \ddot{X} &= \frac{1}{t^\alpha} \ddot{W} - \frac{2\alpha}{t^{\alpha+1}}\dot{W} + \frac{\alpha(\alpha+1)}{t^{\alpha+2}}W.
\end{align*}
Plugging these to \eqref{eq:AGMODEWithr} we get ODE
\begin{align*}
   0 &= \frac{1}{t^\alpha}\ddot{W} 
        + \frac{r-2\alpha}{t^{\alpha+1}} \dot{W} 
        + \frac{\alpha(\alpha +1 - r)}{t^{\alpha+2}} W  
        + \nabla f\pr{ \frac{W}{t^\alpha} + X_\star } .
\end{align*}
Now by defining 
\begin{align*}
    U(W,t) = \frac{\alpha(\alpha +1 - r)}{2t^{\alpha+2}} \norm{W}^2 + t^\alpha \pr{ f \pr{ \frac{W}{t^\alpha} + X_\star} - f_\star }
\end{align*}
we can rewrite the ODE as
\begin{align}          \label{eq:ODE with W for r}
    0 &= \frac{1}{t^\alpha}\ddot{W} +\frac{r-2\alpha}{t^{\alpha+1}}\dot{W} + \nabla_{W} U(W,t).
\end{align}
Now plugging $a(t)=\frac{1}{t^\alpha}$, $b(t)=\frac{r-2\alpha}{t^{\alpha+1}}$,
from conservation law \eqref{eq : Conservation Law for ODE} we get
\begin{align*}   
    E   &\equiv \frac{1}{2t_0^{\alpha}} \norm{\dot{W}(t_0)}^2 
            + \frac{\alpha(\alpha +1 - r)}{2t_0^{\alpha+2}} \norm{W(t_0)}^2 + t_0^\alpha \pr{ f \pr{ \frac{W(t_0)}{t_0^\alpha} + X_\star} - f_\star }     \\
        &= \frac{1}{2 t^\alpha} \norm{\dot{W}}^2 
            + \frac{\alpha(\alpha +1 - r)}{2t^{\alpha+2}} \norm{W}^2 + t^\alpha \pr{ f \pr{ \frac{W}{t^\alpha} + X_\star} - f_\star }
            + \int_{t_{0}}^{t} \frac{2r-3\alpha}{2 s^{\alpha+1}} \norm{\dot{W}}^2 ds \\
        &\quad
        - \int_{t_{0}}^{t} \pr{ \alpha s^{\alpha-1} \pr{  f \pr{ \frac{W}{s^\alpha} + X_\star} - f_\star
                - \inner{ \nabla f \pr{ \frac{W}{s^\alpha} + X_\star } } { \frac{W}{s^\alpha} } } 
                - \frac{\alpha(\alpha +1 - r)(\alpha+2)}{2s^{\alpha+3}} \norm{W}^2 } ds  .
\end{align*}
Rewriting in terms of $X$, $\dot{X}$, $\ddot{X}$ with some reordering we have
\begin{align}   \tag{\ref{eq:EnergyGeneralizedForAlphaR}}
    E   &\equiv  t_0^{\alpha} \pr{ f \pr{ X(t_0) } - f_\star}        \nonumber
        + \frac{1}{2} t_0^{\alpha-2}\norm{ t_0\dot{X}(t_0) + \alpha (X(t_0)-X_\star) }^2   
        + \frac{\alpha(\alpha +1 - r)}{2} t_0^{\alpha-2} \norm{X(t_0)-X_\star}^2   \\
        &= t^{\alpha} \pr{ f \pr{ X } - f_\star}        \nonumber
        + \frac{1}{2} t^{\alpha-2}\norm{ t\dot{X} + \alpha (X-X_\star) }^2   
        + \frac{\alpha(\alpha +1 - r)}{2} t^{\alpha-2} \norm{X-X_\star}^2  \\ &\quad \nonumber
        + \int_{t_0}^{t} \bigg( \frac{(2r-3\alpha)s^{\alpha-3}}{2} \norm{s \dot{X} + \alpha(X-X_\star)}^2 
        + \frac{\alpha(\alpha +1 - r)(\alpha+2)}{2} s^{\alpha-3} \norm{X-X_\star}^2 \bigg) ds \\  &\quad \nonumber
        + \int_{t_0}^{t} \alpha s^{\alpha-1} \pr{  f_\star - f (X)
            - \inner{ \nabla f (X) } { X_\star - X } }  ds    .
\end{align}

\subsection{AGM ODE with $r>3$} \label{appendix : r>3 detail}
Plugging $\alpha=2$, $t_0=0$ to \eqref{eq:EnergyGeneralizedForAlphaR}, we have
\begin{align*}
    E   &\equiv  (5-r) \norm{X_0-X_\star}^2  \\
        &= t^{2} \pr{ f \pr{ X } - f_\star}        
        + \frac{1}{2}\norm{ t\dot{X} + 2 (X-X_\star) }^2   
        + (3-r) \norm{X-X_\star}^2  \\ &\quad 
        + \int_{0}^{t} \bigg( \frac{r-3}{s} \norm{s \dot{X} + 2(X-X_\star)}^2 
                                + \frac{4(3 - r)}{s} \norm{X-X_\star}^2 \bigg) ds \\  &\quad
        + \int_{0}^{t} 2 s \pr{  f_\star - f (X) - \inner{ \nabla f (X) } { X_\star - X } }  ds    .
\end{align*}
Also, since
\begin{align*}
    &\int_{0}^{t} \pr{ 
        \frac{r-3}{s} \norm{s \dot{X} + 2(X-X_\star)}^2 + \frac{4(3 - r)}{s} \norm{X-X_\star}^2 } ds    \\
    &= \int_{0}^{t} \pr{ \frac{r-3}{s} \norm{s \dot{X} }^2 + 4(r-3)\inner{\dot{X}}{X-X_\star} } ds      
    = \int_{0}^{t} \frac{r-3}{s} \norm{s \dot{X} }^2 ds + \br{ 2(r-3)\norm{X-X_\star}^2 }_{0}^{t}      \\
    &= \int_{0}^{t} \frac{r-3}{s} \norm{s \dot{X} }^2 ds + 2(r-3) \pr{ \norm{X-X_\star}^2 - \norm{X_0-X_\star}^2 } .
\end{align*}
Therefore
\begin{align*}
    E   &\equiv  (5-r) \norm{X_0-X_\star}^2  \\
        &= t^{2} \pr{ f \pr{ X } - f_\star}        
        + \frac{1}{2} \norm{ t\dot{X} + 2 (X-X_\star) }^2   
        + (r-3) \norm{X-X_\star}^2  - 2(r-3) \norm{X_0 - X_\star}^2  \\ &\quad 
        + \int_{0}^{t} \frac{r-3}{s} \norm{s \dot{X} }^2 ds 
        + \int_{0}^{t} 2 s \big(  f_\star - f (X) - \inner{ \nabla f (X) } { X_\star - X } \big)  ds    .
\end{align*}

\subsection{AGM ODE with growth condition}      \label{appendix : AGM ODE growth condition}
Rescaling \eqref{eq:ODE with W for r} by multiplying $t^\beta$ we get
\begin{equation*}
    0 = \frac{1}{t^{\alpha-\beta}}\ddot{W} + \frac{r-2\alpha}{t^{\alpha-\beta+1}} \dot{W} + \nabla_{W} \pr{\frac{\alpha(\alpha+1-r)}{2t^{\alpha-\beta+2}} \norm{W}^2
                +  t^{\alpha+\beta} \pr{ f \pr{ \frac{W}{t^\alpha}+X_\star } - f_\star }}.
\end{equation*}
Now plugging $a(t)= \frac{1}{t^{\alpha-\beta}}$, $b(t)=\frac{r-2\alpha}{t^{\alpha-\beta+1}}$,
from conservation law \eqref{eq : Conservation Law for ODE} we get
\begin{align*}
    E   &\equiv \frac{1}{2t_0^{\alpha - \beta}} \norm{ \dot{W}(t_0) }^2 
            + \frac{\alpha(\alpha+1-r)}{2t_0^{\alpha-\beta+2}} \norm{W(t_0)}^2 
            + t_0^{\alpha+\beta} \pr{ f \pr{ \frac{W(t_0)}{t_0^\alpha}+X_\star} - f_\star}    \\
        &= \frac{1}{2t^{\alpha-\beta}} \norm{\dot{W}}^2 
            + \frac{\alpha(\alpha+1-r)}{2t^{\alpha-\beta+2}} \norm{W}^2 
            + t^{\alpha+\beta} \pr{ f \pr{ \frac{W}{t^\alpha}+X_\star} - f_\star }         \\
        &\quad
        + \int_{t_0}^{t} \frac{2r-3\alpha-\beta}{s^{\alpha-\beta+1}} \norm{\dot{W}}^2 ds
        + \int_{t_0}^{t} \frac{\alpha(\alpha+1-r)(\alpha-\beta+2)}{2s^{\alpha-\beta+3}} \norm{W}^2 \, ds \\ &\quad
        - \int_{t_0}^{t} s^{\alpha+\beta-1} \bigg( (\alpha+\beta) \pr{ f \pr{ \frac{W}{s^\alpha}+X_\star} - f_\star }
            - \alpha \inner{ \nabla f \pr{ \frac{W}{s^\alpha}+X_\star} } { \frac{W}{s^\alpha} } \bigg)  ds .
\end{align*}
Rewriting in terms of $X$ we have
\begin{align}       \label{eq:conservation law for rescaled ODE}
    E \nonumber
    &\equiv t_0^{\alpha + \beta}(f(X(t_0))) - f_\star)
        + \frac{1}{2} t_0^{\alpha + \beta - 2} \norm{ t_0\dot{X}(t_0) + \alpha (X(t_0) - X_\star) }^2 
        + \frac{1}{2} \alpha (\alpha + 1 - r) t_0^{\alpha + \beta - 2} \norm{ X(t_0) - X_\star }^2  \\ \nonumber    
    &=  t^{\alpha + \beta}(f(X) - f_\star)
        + \frac{1}{2} t^{\alpha + \beta - 2} \norm{ t\dot{X} + \alpha (X- X_\star) }^2 
        + \frac{1}{2} \alpha (\alpha + 1 - r) t^{\alpha + \beta - 2} \| X - X_\star\|^2     \\  \nonumber   &\quad
        + \int_{t_0}^t \frac{2r - 3\alpha - \beta}{2} s^{\alpha + \beta - 3} \norm{ s\dot{X} + \alpha (X - X_\star) }^2
        + \int_{t_0}^t \frac{\alpha(\alpha + 1 - r)(\alpha - \beta + 2)}{2} s^{\alpha + \beta - 3} \| X - X_\star \|^2 ds \\   &\quad
        + \int_{t_0}^t s^{\alpha + \beta - 1} \bigg( 
            (\alpha + \beta) (f_\star-f(X)) - \alpha \inner{ \nabla f(X) }{ X_\star - X }  \bigg) ds .
\end{align}
To utilize the $\mathbf{H}_1(\gamma)$ hypothesis, it is natural to choose $\alpha, \beta$ such that $\frac{\alpha}{\alpha+\beta} = \frac{1}{\gamma}$.
The choice $\alpha = \frac{2r}{\gamma+2}$, $\beta = \frac{2(\gamma - 1)r}{\gamma + 2}$ makes 
$\frac{\alpha}{\alpha+\beta} = \frac{1}{\gamma}$, and $2r-3\alpha-\beta = 0$, and we get the conservation law used in Section~\ref{sec:growth condition}.
\begin{align*}
    E 
    &\equiv t^{\frac{2\gamma r}{\gamma+2}}(f(X) - f_\star)
        + \frac{1}{2} t^{\frac{2\gamma r}{\gamma+2} - 2} \norm{ t\dot{X} + \frac{2r}{\gamma+2} (X- X_\star) }^2
        + \frac{r( 2 - \gamma (r - 1) )}{(\gamma + 2)^2} t^{\frac{2\gamma r}{\gamma+2} - 2} \| X - X_\star\|^2 \\  &\,\,
     + \int_{t_0}^t \frac{ 2 r ( 2r + 2 - \gamma (r - 1) ) (2 - \gamma (r - 1))}{(\gamma + 2)^3} s^{\frac{2\gamma r}{\gamma+2} - 3} \| X - X_\star \|^2 ds \\ &\,\,
    +  \int_{t_0}^t s^{\frac{2\gamma r}{\gamma+2} - 1} \frac{2\gamma r}{\gamma+2}
            \left(  f_\star - f(X) - \frac{1}{\gamma} \inner{ \nabla f(X) }{ X_\star - X } \right) ds .
\end{align*}

\subsubsection{Lyapunov function for $r>3$ in \cite{SuBoydCandes2014_differential}}
Plugging $\alpha=r-1$, $\beta=3-r$, $t_0=0$ to \eqref{eq:conservation law for rescaled ODE}, 
we have
\begin{align*}
    E 
    &\equiv \frac{(r-1)^2}{2} \norm{ X_0 - X_\star }^2  \\     
    &=  t^{2}(f(X) - f_\star)
        + \frac{1}{2} \norm{ t\dot{X} + (r-1) (X- X_\star) }^2      \\     &\quad
        + \int_{0}^t s (r-1) \bigg( 
            f_\star-f(X) - \inner{ \nabla f(X) }{ X_\star - X }  \bigg) ds 
        + \int_{0}^t s (r-3) (f(X) - f_\star) ds  .
\end{align*}
Since all terms are nonnegative, we immediately get
\begin{align*}
    f(X) - f_\star \le \frac{(r-1)^2}{2t^2} \norm{ X_0 - X_\star }^2.
\end{align*}
In \cite{SuBoydCandes2014_differential}, they also present
\begin{align*}
    \int_{0}^{\infty} t(f(X(t))-f^\star) dt \le \frac{(r-1)^2}{2(r-3)} \norm{X_0-X^\star}^2,
\end{align*}
and this can also be obtained immediately from conservation law.

\subsection{SC-AGM ODE} \label{appendix : SC-AGM detail}
We proceed the argument similar to \ref{appendix : derivation of mechanical energy for general r}.
Start with the ODE \eqref{eq:AGMODEWithr}
\begin{align*}
    0 = \ddot{X} + 2\sqrt{\mu} \dot{X} + \nabla f(X).
\end{align*}
Now consider the coordinate change $W=e^{\beta t}(X-X_\star)$.
Then we see
\begin{align*}
    W           &= e^{\beta t} (X-X_\star)    \\
    \dot{W}     &= e^{\beta t} \pr{ \dot{X} + \beta (X-X_\star) }    \\
    \ddot{W}    &= e^{\beta t} \pr{ \ddot{X} + 2\beta \dot{X} + \beta^2 (X-X_\star) } .    
\end{align*}
From this, we can rewrite $X$, $\dot{X}$, $\ddot{X}$ in terms of $W$, $\dot{W}$, $\ddot{W}$,
\begin{align*}
    X           &= e^{-\beta t} W + X_\star  \\
    \dot{X}     &= e^{-\beta t} \pr{ \dot{W} - \beta W}  \\
    \ddot{X}    &= e^{-\beta t} \pr{ \ddot{W} - 2\beta \dot{W} + \beta^2 W } .
\end{align*}
Plugging these to \eqref{eq:AGMODEWithr} we get ODE
\begin{align*}
   0 &= e^{-\beta t} \pr{ \ddot{W} + 2(\sqrt{\mu}-\beta) \dot{W} + \beta(\beta-2\sqrt{\mu}) W } 
        + \nabla f\pr{e^{-\beta t} W + X_\star} .
\end{align*}
Now by defining 
\begin{align*}
    U(W,t) = \frac{\beta(\beta-2\sqrt{\mu})}{2} e^{-\beta t} \norm{W}^2 + e^{\beta t} \pr{ f\pr{e^{-\beta t} W + X_\star} - f_\star},
\end{align*}
we can rewrite the ODE as
\begin{align*}
    0 = e^{-\beta t} \ddot{W} + 2(\sqrt{\mu}-\beta) e^{-\beta t} \dot{W} + \nabla_{W} U(W,t).
\end{align*}
Now plugging $a(t)=e^{-\beta t}$, $b(t)=2(\sqrt{\mu}-\beta) e^{-\beta t}$,
from conservation law \eqref{eq : Conservation Law for ODE} we get
\begin{align*}   
    E   &\equiv \frac{e^{-\beta t_0}}{2} \norm{\dot{W}(t_0)}^2 
            + \frac{\beta(\beta-2\sqrt{\mu})}{2} e^{-\beta t_0} \norm{W(t_0)}^2 + e^{\beta t_0} \pr{ f(e^{-\beta t_0}W(t_0) + X_\star) - f_\star}     \\
        &= \frac{e^{-\beta t}}{2} \norm{\dot{W}}^2 
            + \frac{\beta(\beta-2\sqrt{\mu})}{2} e^{-\beta t} \norm{W}^2 + e^{\beta t} \pr{ f\pr{e^{-\beta t} W + X_\star} - f_\star}
            + \int_{t_{0}}^{t} \frac{4\sqrt{\mu} - 3\beta}{2}  e^{-\beta s} \norm{\dot{W}}^2 ds \\
        &\quad
        - \int_{t_{0}}^{t} \pr{ \beta e^{\beta s} \pr{  f(e^{-\beta s}W + X_\star) - f_\star
                - \inner{ \nabla f(e^{-\beta s}W + X_\star) } { e^{-\beta s}W } } 
                - \frac{\beta^2(\beta-2\sqrt{\mu})}{2} e^{-\beta s} \norm{W}^2 } ds.
\end{align*}
Plugging $t_0=0$ and rewriting in terms of $X$, $\dot{X}$, $\ddot{X}$ we have
\begin{align*}   
    E   &\equiv  f \pr{ X_0 } - f_\star + \beta(\beta-\sqrt{\mu}) \norm{X_0-X_\star}^2    \\
        &= e^{\beta t} \pr{ f \pr{ X } - f_\star       
        + \frac{1}{2} \norm{\dot{X} + \beta(X-X_\star)}^2   
        + \frac{\beta(\beta-2\sqrt{\mu})}{2} \norm{X-X_\star}^2 } \\ &\quad 
        + \int_{0}^{t} \frac{4\sqrt{\mu} - 3\beta}{2} e^{\beta s} \norm{\dot{X} + \beta(X-X_\star)}^2 ds \\&\quad
        + \int_{0}^{t} \beta e^{\beta s} 
            \pr{  f_\star - f (X)  - \inner{ \nabla f (X) } { X_\star - X } 
            + \frac{\beta(\beta-2\sqrt{\mu})}{2} \norm{X-X_\star}^2 }  ds    .
\end{align*}
Now plugging $\beta=\sqrt{\mu}$ we have
\begin{align*}   
    E   &\equiv  f \pr{ X_0 } - f_\star     \\
        &= e^{\sqrt{\mu} t} \pr{ f \pr{ X } - f_\star       
        + \frac{1}{2} \norm{\dot{X} + \sqrt{\mu}(X-X_\star)}^2   
        - \frac{\mu}{2} \norm{X-X_\star}^2 } \\ &\quad 
        + \int_{0}^{t} \frac{\sqrt{\mu}}{2} e^{\sqrt{\mu} s} \norm{\dot{X} + \sqrt{\mu}(X-X_\star)}^2 ds \\  &\quad
        + \int_{0}^{t} \sqrt{\mu} e^{\sqrt{\mu} s} 
            \pr{  f_\star - f (X) - \inner{ \nabla f (X) } { X_\star - X } 
            - \frac{\mu}{2} \norm{X-X_\star}^2 }  ds    .
\end{align*}
Finally, from
\begin{align*}
    & \int_{0}^{t} \frac{\sqrt{\mu}}{2} e^{\sqrt{\mu} s} \norm{\dot{X} + \sqrt{\mu}(X-X_\star)}^2 ds    \\ &\quad
    = \int_{0}^{t} \pr{ \frac{\sqrt{\mu}}{2} e^{\sqrt{\mu} s} \norm{\dot{X}}^2 
      + \frac{\mu}{2} e^{\sqrt{\mu} s} \pr{ 2\inner{\dot{X}}{X-X_\star} + \sqrt{\mu} \norm{X-X_\star}^2 } } ds      \\&\quad
    = \int_{0}^{t} \pr{ \frac{\sqrt{\mu}}{2} e^{\sqrt{\mu} s} \norm{\dot{X}}^2 
      + \frac{\mu}{2} \frac{d}{ds} \pr{ e^{\sqrt{\mu} s} \norm{X-X_\star}^2 } } ds      \\&\quad
    = \int_{0}^{t} \frac{\sqrt{\mu}}{2} e^{\sqrt{\mu} s} \norm{\dot{X}}^2  ds  
        +  \frac{\mu}{2} \br{ e^{\sqrt{\mu} s} \norm{X-X_\star}^2 }_{0}^{t}  \\&\quad
    = \int_{0}^{t} \frac{\sqrt{\mu}}{2} e^{\sqrt{\mu} s} \norm{\dot{X}}^2  ds  
        +  \frac{\mu}{2} \pr{ e^{\sqrt{\mu} t} \norm{X-X_\star}^2 - \norm{X_0-X_\star}^2 }.
\end{align*}
we conclude
\begin{align*}   
    E   &\equiv  f \pr{ X_0 } - f_\star     \\
        &= e^{\sqrt{\mu} t} \pr{ f \pr{ X } - f_\star       
        + \frac{1}{2} \norm{\dot{X} + \sqrt{\mu}(X-X_\star)}^2    }
        - \frac{\mu}{2} \norm{X_0-X_\star}^2 \\ &\quad 
        + \int_{0}^{t} \frac{\sqrt{\mu}}{2} e^{\sqrt{\mu} s} \norm{\dot{X}}^2 ds 
        + \int_{0}^{t} \sqrt{\mu} e^{\sqrt{\mu} s} 
            \pr{  f_\star - f (X) - \inner{ \nabla f (X) } { X_\star - X } 
            - \frac{\mu}{2} \norm{X-X_\star}^2 }  ds    .
\end{align*}

\subsection{Gradient flow}
Recall, gradient flow was written as
\begin{align*}
    0 = \dot{X} + \nabla f(X).
\end{align*}
Consider the dilated coordinate $W=t(X-X_\star)$. Then we see
\begin{align*}
    W       &= t(X-X_\star)   \\
    \dot{W} &= t\dot{X} + (X-X_\star)   .
\end{align*}
Then $X$, $\dot{X}$ can be rewritten as
\begin{align*}
    X       &= \frac{W}{t} + X_\star    \\
    \dot{X} &= \frac{\dot{W}}{t} - \frac{W}{t^2}   .
\end{align*}
Plugging these to ODE, we have
\begin{align*}
    0 = \frac{\dot{W}}{t} - \frac{W}{t^2} + \nabla f \pr{ \frac{W}{t} + X_\star } .
\end{align*}
Now by defining
\begin{align*}
    U(W,t) = -\frac{1}{2t^2} \norm{W}^2 + t \pr{ f \pr{ \frac{W}{t} + X_\star } - f_\star }  ,
\end{align*}
we can rewrite ODE as
\begin{align*}
    0 = \frac{\dot{W}}{t} + \nabla_W U(W,t) .
\end{align*}
Now plugging $a(t)=0$, $b(t)=\frac{1}{t}$, from conservation law \eqref{eq : Conservation Law for ODE}
\begin{align*}
    E   &\equiv \lim_{t_0\to0} U(W(t_0),t_0)  \\
        &= \int_{0}^t \frac{1}{s} \norm{\dot{W}}^2 ds + U(W,t) - \int_{0}^t \pd{s} U(W,s) ds    \\
        &= \int_{0}^t \frac{1}{s} \norm{\dot{W}}^2 ds 
            - \frac{1}{2t^2} \norm{W}^2 + t \pr{ f \pr{ \frac{W}{t} + X_\star } - f_\star }     \\&\quad
            - \int_{0}^t \pr{ \frac{1}{s^3} \norm{W}^2 
                + \pr{  f \pr{ \frac{W}{s} + X_\star } - f_\star
                + s \inner{ \nabla f \pr{ \frac{W}{s} + X_\star } } { -\frac{W}{s^2} } }} ds  .
\end{align*}
Rewriting in terms of $X$, $\dot{X}$, we get the conservation law in Section~\ref{s:gradient flow}
\begin{align*}
    E   &\equiv -\frac{1}{2}\norm{X_0-X_\star}^2 \\
        &= t \pr{ f (X) - f_\star } - \frac{1}{2} \norm{X-X_\star}^2        \\&\quad
            +  \int_{0}^t \pr{ \frac{1}{s} \norm{ s\dot{X} + (X-X_\star)}^2 - \frac{1}{s} \norm{X-X_\star}^2 }ds
            -  \int_{0}^t \pr{ f(X) - f_\star - \inner{\nabla f(X)}{X-X_\star} }ds    \\
        &= t \pr{ f (X) - f_\star } - \frac{1}{2} \norm{X-X_\star}^2        \\&\quad
            +  \int_{0}^t \pr{ s \norm{ \dot{X} } ^2 + \frac{d}{ds}\norm{X-X_\star}^2 }ds
            +  \int_{0}^t \pr{ f_\star - f(X) - \inner{\nabla f(X)}{X_\star-X} }ds    \\
        &= t \pr{ f (X) - f_\star } + \frac{1}{2} \norm{X-X_\star}^2 - \norm{X_0-X_\star}^2
            +  \int_{0}^t s \norm{ \dot{X} } ^2 ds
            +  \int_{0}^t \pr{ f_\star - f(X) - \inner{\nabla f(X)}{X_\star-X} }ds    .
\end{align*}

\section{Omitted calculations of Section~\ref{s:section4}}
\label{s:sectionC}

\subsection{Derivation of OGM-G ODE} \label{appendix : OGM-G ODE derivation}
OGM-G in \cite{KimFessler2021_optimizing} was presented as 
\begin{align*}
    x_k^+   &= x_k - \frac{1}{L}\nabla f(x_k)       \\
    x_{k+1} &= x_k^+    + \frac{(\theta_{K-k}-1)(2\theta_{K-(k+1)}-1)}{\theta_{K-k}(2\theta_{K-k}-1)} ( x_k^+ - x_{k-1}^+ ) 
                        + \frac{2\theta_{K-(k+1)}-1}{2\theta_{K-k}-1} ( x_k^+ - x_k ) .
\end{align*}
Plugging $ x_k^+ = x_k - \frac{1}{L}\nabla f(x_k) $ to the second line and
using the fact $\theta_{K-k} = \frac{K-k}{2} + o(K)$ we have
\begin{align*}
    x_{k+1} 
        &= x_k - \frac{1}{L}\nabla f(x_k)   
            + \frac{(K-k-2+ o(K))^2}{(K-k+ o(K)) (K-k-1+ o(K)) } \pr{ x_k - x_{k-1} - \frac{1}{L} \pr{ \nabla f(x_k) - \nabla f(x_{k-1}) } }   \\&\quad
            - \frac{K-k-2+o(K)}{K-k-1+o(K)} \frac{1}{L}\nabla f(x_k)  \\
        &= x_k    
            + \pr{ 1 - \frac{3(K-k)+o(K)}{(K-k)^2+o(K)K} } \pr{ x_k - x_{k-1} }   
            - \pr{ 2 - \frac{1}{K-k+o(K)} } \frac{1}{L}\nabla f(x_k)  \\&\quad
            - \frac{1}{L} \frac{(K-k)^2+ o(K)K}{(K-k)^2+ o(K)K} (\nabla f(x_k) - \nabla f(x_{k-1})).
\end{align*}
Similar to \cite{SuBoydCandes2014_differential}, 
we use the identification $\frac{1}{L} = h^2$, $t=kh$ and $x_k=X(kh)$.
Moreover for fixed $T>0$, we use identification $T=Kh$.
Adding $-2x_k + x_{k-1}$ and dividing $h^2$ both sides we have
\begin{align*}
    \frac{ (x_{k+1} - x_k) - (x_k-x_{k-1}) }{ h^2 }   
        &=  -   \frac{3(Kh-kh)+o(K)h}{(Kh-kh)^2+o(K)Kh^2} \frac{x_k - x_{k-1}}{h}     
            -   \pr{ 2 - \frac{h}{Kh-kh+o(K)h} } \nabla f(x_k)      \\
        &\quad    -  \frac{(Kh-kh)^2+ o(K)Kh^2}{(Kh-kh)^2+ o(K)Kh^2} (\nabla f(x_k) - \nabla f(x_{k-1})) \\
        &=  -   \frac{3}{T-t} \frac{X(t) - X(t-h)}{h}     
            -   2 \nabla f(X(t))         
            -  (\nabla f(X(t)) - \nabla f(X(t-h)) 
            + o(K)h.
\end{align*}
Finally taking limit $h \to 0$, we obtain the desired ODE 
\begin{align*}
    0 = \ddot{X}(t) - \frac{3}{t-T}\dot{X}(t) + 2\nabla f(X(t)).
\end{align*}

\subsubsection{OGM-G ODE coincides with the ODE model of OBL-G${}_\flat$}
The method OBL-G${}_\flat$ \cite{ParkRyu2021_optimal} 
\begin{align*}
    x_{k}^{+} &= x_k - \frac{1}{L} \nabla f(x_k)  \\
    z_{k+1} &= z_k - \frac{1}{L} \frac{K-k+1}{2} \nabla f(x_k)  \\
    x_{k+1} &= \frac{K-k-2}{K-k+2} x_{k}^{+} + \frac{4}{K-k+2} z_{k+1} .
\end{align*}
is a variant of OGM-G.
Interestingly, the ODE model of OBL-G${}_\flat$ exactly coincides with OGM-G ODE.

Note this method is written in the form with auxiliary sequence $z_k$, we derive the ODE in a different way.
We take the same identification $\frac{1}{L} = h^2$, $Kh=T$, $kh=t$, $x_k=X(kh)$, $z_k=Z(kh)$. 
Then we may regard the method as a system of first-order ODEs.
From $z_k$ update, by taking limit $h \to 0$ we have
\begin{align*}   
    \frac{z_{k+1}-z_k}{h} = -\frac{Kh-kh+h}{2} \nabla f(x_k)  
    \quad \stackrel{h \to 0}{\Longrightarrow}
    \quad \dot{Z}(t) = -\frac{T-t}{2} \nabla f(X) .
\end{align*}
From $x_k$ update, dividing both sides by $h$, subtracting $x_{k}^{+}$ both sides and by taking limit $h \to 0$
we have
\begin{align}   \label{eq : OBL-G dotZ}
    \frac{x_{k+1}-x_k}{h}  
    &= \frac{4}{Kh-kh+2h} ( z_{k+1} - x_k ) - \frac{Kh-kh-2h}{Kh-kh+2h} \nabla f(x_k) h 
    \quad \stackrel{h \to 0}{\Longrightarrow}
    \quad \dot{X}(t) = \frac{4}{T-t} (Z(t)-X(t)) .
\end{align}
Thus we get system of first-order ODEs.
Now to derive a second-order ODE, multiplying $T-t$ to \eqref{eq : OBL-G dotZ} and differentiating, we have
\begin{align*}
    (T-t)\ddot{X}(t) - \dot{X}(t)
    = 4\pr{ \dot{Z}(t) - \dot{X}(t) }
    = 4\pr{ -\frac{T-t}{2} \nabla f(X) - \dot{X}(t) }.
\end{align*}
Dividing $T-t$ and organizing the result, we conclude
\begin{align*}
    0 = \ddot{X}(t) - \frac{3}{t-T} \dot{X}(t) + 2 \nabla f(X).
\end{align*}

\subsection{Conservation law for OGM-G ODE} \label{appendix : Energy for OGM-G ODE }
We proceed argument similar to 
\ref{appendix : SC-AGM detail}.
Start with ODE presented in Section~\ref{section : OGM-G r<-3}   
\begin{align}      \label{eq:generalized_OGMG_ODE}
    0 = \ddot{X} + \frac{r}{t-T} \dot{X} + 2\nabla f(X).
\end{align}
Now consider the coordinate change $W=(T-t)^{\alpha}(X-X_c)$.

Then we see
\begin{align*}
    W(t)        &= 
        (T-t)^\alpha (X(t) - X_c)    \\
    \dot{W}(t)  &= 
        (T-t)^\alpha \dot{X}(t) - \alpha (T-t)^{\alpha-1} (X(t) - X_c)   \\
    \ddot{W}(t) &= 
        (T-t)^\alpha \ddot{X}(t) - 2 \alpha (T-t)^{\alpha-1} \dot{X}(t) + \alpha(\alpha-1) (T-t)^{\alpha-2} (X(t) - X_c).
\end{align*}
Note the sign flips while differentiating $(T-t)^\alpha$.

From this, we can rewrite $X$, $\dot{X}$, $\ddot{X}$ in terms of $W$, $\dot{W}$, $\ddot{W}$,
\begin{align*}
    X(t)        
        &= (T-t)^{-\alpha} W(t)   + X_c  \\
    \dot{X}(t)  
        &= (T-t)^{-\alpha} \dot{W}(t) + \alpha (T-t)^{-\alpha-1} W(t)   \\
    \ddot{X}(t) 
        &= (T-t)^{-\alpha} \ddot{W}(t)  + 2\alpha (T-t)^{-\alpha-1} \dot{W}(t)
            + \alpha(\alpha+1) (T-t)^{-\alpha-2} W(t) .
\end{align*}
Plugging these to \eqref{eq:AGMODEWithr} we get ODE
\begin{align*}
    0   &= \frac{1}{(T-t)^{\alpha}}\ddot{W} 
            + \frac{2\alpha-r}{(T-t)^{\alpha+1}} \dot{W} 
            + \frac{\alpha(\alpha + 1 - r )}{(T-t)^{\alpha+2}} W  
            + 2\nabla f \pr{ \frac{W}{(T-t)^\alpha}+X_c}   .
\end{align*}
Now by defining 
\begin{align*}
    U(W,t) = \frac{\alpha(\alpha+1-r)}{2(T-t)^{\alpha+2}} \norm{W}^2
            + 2(T-t)^{\alpha} \pr{ f \pr{ \frac{W}{(T-t)^\alpha}+X_c} - f(X_c) }
\end{align*}
we can rewrite the ODE as
\begin{align*}
    0 = \frac{1}{(T-t)^{\alpha}}\ddot{W} 
            + \frac{2\alpha-r}{(T-t)^{\alpha+1}} \dot{W} 
            + \nabla_{W} U(W,t) .
\end{align*}
Now plugging $a(t)=\frac{1}{(T-t)^{\alpha}}$, $b(t)=\frac{2\alpha-r}{(T-t)^{\alpha+1}}$,
from conservation law \eqref{eq : Conservation Law for ODE} we get
\begin{align*}   
    E   &\equiv \frac{1}{2(T-t_0)^{\alpha}} \norm{\dot{W}(t_0)}^2 
            + \frac{\alpha(\alpha+1-r)}{2(T-t_0)^{\alpha+2}} \norm{W(t_0)}^2 
            + 2(T-t_0)^{\alpha} \pr{ f \pr{ \frac{W(t_0)}{(T-t_0)^\alpha}+X_c} - f(X_c) }     \\
        &= \frac{1}{2(T-t)^{\alpha}} \norm{\dot{W}}^2 
            + \frac{\alpha(\alpha+1-r)}{2(T-t)^{\alpha+2}} \norm{W}^2 
            + 2(T-t)^{\alpha} \pr{ f \pr{ \frac{W}{(T-t)^\alpha}+X_c} - f(X_c) }         \\
        &\quad
        + \int_{t_{0}}^{t} \frac{3\alpha-2r}{2(T-s)^{\alpha+3}} \norm{\dot{W}}^2 ds
        - \int_{t_0}^{t} \frac{\alpha(\alpha +1 - r)(\alpha+2)}{2(T-s)^{\alpha+3}}  \norm{W}^2 \, ds \\ &\quad
        - \int_{t_0}^{t} \frac{2\alpha}{(T-s)^{\alpha+1}} \pr{  f(X_c) - f \pr{ \frac{W}{(T-s)^\alpha}+X_c}  - \inner{ \nabla f \pr{ \frac{W}{(T-s)^\alpha}+X_c} } { \frac{W}{(T-s)^\alpha} } }  ds .
\end{align*}
Plugging $t_0=0$ and rewriting in terms of $X$, $\dot{X}$, $\ddot{X}$ we have
\begin{align}   \label{eq : Energy for OGM-G for general r, alpha}
    E   &= 2T^{\alpha} \pr{ f (X_0) - f(X_c)}        
        + \pr{ \frac{\alpha^2}{2} + \frac{\alpha(\alpha +1 - r)}{2}} T^{\alpha-2} \norm{ X_0-X_c }^2     \\   \nonumber
        &= 2(T-t)^{\alpha} \pr{ f \pr{ X } - f(X_c)}        
        + \frac{1}{2} (T-t)^{\alpha-2}\norm{ (T-t)\dot{X} - \alpha (X-X_c) }^2   
        + \frac{\alpha(\alpha +1 - r)}{2} (T-t)^{\alpha-2} \norm{X-X_c}^2  \\  \nonumber &\quad
        + \int_{0}^{t} \pr{ \frac{3\alpha-2r}{2} (T-s)^{\alpha-3} \norm{(T-s) \dot{X} - \alpha(X-X_c)}^2    
                - \frac{\alpha(\alpha +1 - r)(\alpha+2)}{2} (T-s)^{\alpha-3} \norm{X-X_c}^2 } \; ds  \\ &\quad
        + \int_{0}^{t} (-2\alpha) (T-s)^{\alpha-1} \pr{  f(X_c) - f (X)  - \inner{ \nabla f (X) } { X_c - X } }  ds \nonumber .
\end{align}
Now plugging $\alpha=-2$ we get the energy in Section 4.2, 
moreover plugging $r=-3$ we get the energy for $r=-3$ in Section 4.

\subsection{Regularity of OGM-G ODE at terminal time $T$} \label{appendix : Property of OGM-G ODE}
Since the argument for $r=-3$ is exactly same for general $r$, 
we prove the statement for the general $r<0$.
We will present our proofs in following order.
\begin{itemize}
    \item[(i)]      $ \sup_{t\in[0,T)}\norm{\dot{X}(t)}$ is bounded.
    \item[(ii)]     $ X(t)$ can be continuously extended to $T$.
    \item[(iii)]    $ \lim_{t\to T^{-}} \dot{X}(t) = 0$.
    \item[(iv)]     $ \lim_{t\to T^{-}} \frac{\dot{X}(t)}{t-T} = -\frac{2}{1+r} \nabla f(X(T))$.
    \item[(v)]      $ \lim_{t\to T^{-}} \ddot{X}(t) = -\frac{2}{1+r} \nabla f(X(T))$.
\end{itemize}    
(i), (ii) holds for $r\le0$, (iii) holds for $r<0$, and (iv), (v) holds for $r<0$ with $r\ne-1$.

\subsubsection{$\sup_{t\in[0,T)}\norm{\dot{X}(t)}$ is bounded if $r\le0$} \label{section : dotX is bounded}
Considering conservation law \eqref{eq : Energy for OGM-G for general r, alpha} with $\alpha=0$, $X_c=X_0$, 
we have
\begin{align}       \label{eq:energy for zero alpha}
    E   &\equiv 0 
        =  \frac{1}{2} \norm{\dot{X}(t)}^2 + 2( f(X(t)) - f(X_0) ) 
            - \int_{0}^{t} \frac{r}{T-s} \norm{\dot{X}(s)}^2 ds .
\end{align}
Collecting the terms except the integrand, define $\Psi \colon [0,T) \to \reals$ as
\begin{align*}
    \Psi(t)=\frac{1}{2} \norm{\dot{X}(t)}^2 + 2( f(X(t)) - f(X_0) ) .
\end{align*}
Observe for $r\le0$ 
    $$ \dot{\Psi}(t) = \frac{r}{T-t} \norm{\dot{X}(t)}^2 \le 0 , $$
so $\Psi(t)$ is a nonincreasing function.
Thus $\Psi(t)\le \Psi(0) = 0$, 
and from the fact $f_\star=\inf_{x\in\reals^n} f(x) > -\infty$,
we have
\begin{align*}
    \norm{\dot{X}(t)}^2
        = 2\Psi(t) + 4(f(X_0)-f(X(t)))
        \le 4(f(X_0)-f_\star) .
\end{align*}
Therefore $\sup_{t\in[0,T)}\norm{\dot{X}(t)} \le 2\sqrt{f(X_0)-f_\star}$, we get the desired result.

\subsubsection{$X(t)$ can be continuously extended to $T$} \label{section : X can be extended continuously}
We first prove $X(t)$ is uniformly continuous.
From the result of \ref{section : dotX is bounded}, we see
\begin{align*}
    \norm{ X(t) - X(t+\delta) }
        = \norm{ \int_{t}^{t+\delta} \dot{X}(s) ds }
        \le \int_{t}^{t+\delta} \norm{ \dot{X}(s) } ds
        \le \int_{t}^{t+\delta} 2\sqrt{f(X_0)-f_\star} ds
        = 2 \delta \sqrt{f(X_0)-f_\star} .
\end{align*}
Thus for $X$ is $2 \sqrt{f(X_0)-f_\star}$-Lipschitz function, 
we can conclude X is uniformly continuous.

Now from the fact of basic analysis,
we know for $D\subset \reals^n$, uniformly continuous function $g:D\to\reals^n$ can be extended continuously to $\bar{D}$.
Therefore $X\colon[0,T)\to\reals^n$ can be extended to $\overline{[0,T)}=[0,T]$, we get the desired result.

\subsubsection{$\lim_{t\to T^{-}}\norm{\dot{X}(t)}=0$}      \label{appendix:dotX_becomes_zero}
We first prove the limit $\lim_{t\to T^{-}}\norm{\dot{X}(t)}$ exists.
From $\Psi$ defined in \ref{section : dotX is bounded} we have 
\begin{align*}
    \norm{\dot{X}(t)} = \sqrt{ 2\Psi(t) + 4(f(X_0)-f(X(t))) },
\end{align*}
so it is enough to show $\lim_{t\to T^{-}}\Psi(t)$ and $\lim_{t\to T^{-}}f(X(t))$ exists.
From \ref{section : X can be extended continuously} we know $\lim_{t\to T^{-}} X(t)$ exists,
thus from continuity of $f$, we have $\lim_{t\to T^{-}}f(X(t))$ exists.
It remains to show $\lim_{t\to T^{-}}\Psi(t)$ exists.

Recall $\Psi$ is nonincreasing. 
Moreover, since $f_\star=\inf_{x\in\reals^n} f(x) > -\infty$ we have
\begin{align*}
    \Psi(t) = \frac{1}{2}\norm{\dot{X}(t)}^2 + 2(f(X(t))-f(X_0))
            \ge 2(f_\star - f(X_0)) ,
\end{align*}
so $\Psi$ is bounded below.
Thus $\Psi$ is nonincreasing and bounded below, by completeness of real numbers,
we conclude $\lim_{t\to T^{-}}\Psi(t)$ exists.
Therefore $\lim_{t\to T^{-}}\norm{\dot{X}(t)}$ exists.

Now we prove $\lim_{t\to T^{-}}\norm{\dot{X}(t)}=0$. 
Let $C=\lim_{t\to T^{-}}\norm{\dot{X}(t)}\ge0$. 
Assume for contradiction that $C>0$. 
Then there is $\epsilon>0$ such that $T-\epsilon < s < T$ implies $\norm{\dot{X}(s)} > \frac{C}{2}$.
Thus for $t >T-\epsilon$, if $r\le0$ we have 
\begin{align*}
    \int_{0}^{t} \frac{r}{T-s} \norm{\dot{X}(s)}^2 ds 
    = \int_{0}^{T-\epsilon} \frac{r}{T-s} \norm{\dot{X}(s)}^2 ds 
        + \int_{T-\epsilon}^{t} \frac{r}{T-s} \norm{\dot{X}(s)}^2 ds
    \le \int_{T-\epsilon}^{t} \frac{C^2}{4} \frac{r}{T-s} ds.
\end{align*} 
Since $\lim_{t\to T^{-}} \int_{T-\epsilon}^{t} \frac{C^2}{4} \frac{r}{(T-s)} ds = -\infty$ if $r<0$,
we conclude $\lim_{t\to T^{-}} \int_{0}^{t} \frac{r}{T-s} \norm{\dot{X}(s)}^2 ds = -\infty$ from above inequality.
By the way from \eqref{eq:energy for zero alpha} 
we know $\Psi(t) = \int_{0}^{t} \frac{r}{T-s} \norm{\dot{X}(s)}^2 ds$, 
but we have just observed above that $\Psi(t)$ is bounded below.
This is a contradiction, we conclude $\lim_{t\to T^{-}}\norm{\dot{X}(t)}=0$.

\subsubsection{$ \lim_{t\to T^{-}} \frac{\dot{X}(t)}{t-T} = -\frac{2}{r+1}\nabla f(X(T))$ } \label{appendix:dotX_t-T_limit}
The key observation of the proof is
\begin{align*}
    \frac{d}{dt} \pr{(T-t)^r\dot{X}(t)} = - 2(T-t)^r \nabla f(X(t)).
\end{align*}
We can check above is true from the ODE $0 = \ddot{X} + \frac{r}{t-T} \dot{X} + 2\nabla f(X)$.
With this observation, we can handle the separated terms $\ddot{X}$ and $\dot{X}$ as one term.

Integrating both sides from $0$ to $t$, we get
\begin{align*}
    (T-t)^r\dot{X}(t) = - \int_{0}^t 2(T-s)^r \nabla f(X(s)) ds.
\end{align*}
Multiplying $(T-t)^{-(r+1)}$, we get
\begin{align}   \label{eq : average rate of change for dotX}
    \frac{\dot{X}(t)}{T-t} = -(T-t)^{-(r+1)} \int_{0}^t 2 (T-t)^r \nabla f(X(s)) ds  .
\end{align}
From \citep[Corollary~25.5.1]{Rockafellar1970_convex}, 
the fact $f$ is convex and differentiable implies continuity of $\nabla f$.
From \ref{section : X can be extended continuously}, we see $\lim_{t\to T^{-}} \nabla f(X(t))$ exists.
Moreover from \ref{appendix:dotX_becomes_zero}, we see the numerator for left hand side reaches to zero as $t \to T^{-}$. 
Therefore we can apply L'H\^{o}pital's rule (componentwisely), for $r\ne-1$ we conclude
\begin{align*}
    \lim_{t\to T^{-}} \frac{\dot{X}(t)}{T-t}
        &= -\lim_{t\to T^{-}} \frac{\int_{0}^t 2 (T-t)^r \nabla f(X(s)) \;ds }{ (T-t)^{r+1} } 
        = \frac{2}{r+1} \lim_{t\to T^{-}} \nabla f(X(t))
        = \frac{2}{r+1} \nabla f(X(T)).
\end{align*}
By flipping the sign of both sides, we get the desired result.

\subsubsection{$ \lim_{t\to T^{-}} \ddot{X}(t) = -\frac{2}{r+1} \nabla f(X(T))$}
From ODE \eqref{eq:generalized_OGMG_ODE} we have
    $$ \ddot{X}(t) = \frac{r}{T-t} \dot{X}(t) - 2 \nabla f(X(t)). $$
We know the limit $t\to T^{-}$ for right hand side exists by \ref{appendix:dotX_t-T_limit}.
Therefore $\lim_{t\to T^{-}}\ddot{X}(t)$ exists, by L'H\^{o}pital's rule we have
\begin{align*}
    \lim_{t\to T^{-}} \ddot{X}(t) 
        = \lim_{t\to T^{-}} \frac{\dot{X}(t)}{t-T}
        = -\frac{2}{r+1} \nabla f(X(T)).
\end{align*}

\subsection{Correspondence with discrete analysis of OGM-G} \label{appendix : OGM-G discrete correpondence}
\citet{LeeParkRyu2021_geometric} presented Lyapunov function proof for convergence analysis of OGM-G. 
They first rewrote OGM-G with auxiliary sequence $z_k$ as follows
\begin{align}
    x_{k}^{+}   &= x_{k} - \frac{1}{L} \nabla f(x_{k})      \nonumber \\
    z_{k+1}     &= z_k - \frac{\theta_{K-k}}{L} \nabla f(x_k)   \label{eq : OGM-G z update}    \\
    x_{k+1}     &= \frac{\theta_{K-(k+2)}^4}{\theta_{K-(k+1)}^4} x_k^{+} \label{eq : OGM-G x update}
                    + \pr{ 1 - \frac{\theta_{K-(k+2)}^4}{\theta_{K-(k+1)}^4} } z_{k+1}.
\end{align}
Then they presented the Lyapunov function as follows
\begin{align}   \label{eq:Discrete OGM-G Lyapunov}
    U_k &= \frac{1}{\theta_{K-k}^2} 
            \pr{ \frac{1}{2L} \norm{f(x_K)}^2 + \frac{1}{2L} \norm{f(x_k)}^2 
            + f(x_k) - f(x_K) - \inner{\nabla f(x_k)}{x_k-x_{k-1}^+} }  \\&\quad
            + \frac{L}{\theta_{K-k}^4} \inner{z_k-x_{k-1}^+}{z_k-x_K^+}. \nonumber
\end{align}
We claim there is a correspondence between this function and the Lyapunov function we've presented in Theorem~\ref{thm:main_result_for_OGMG}.
We use same identification as did in \ref{appendix : OGM-G ODE derivation},
$\frac{1}{L} = h^2$, $kh=t$, $Kh=T$, $x_k=X(kh)$, $z_k=Z(kh)$.
Then we derive continuous counterpart of $U_k$ by dividing $2h^2$ then ignoring $o(K)h$ and $O(h)$.

We first calculate the continuous counterpart of $z_k$.
Rewrite the update equation \eqref{eq : OGM-G z update} as
\begin{align}       \label{eq : rewritten OGM-G x update}
    x_{k+1} - x_k^{+}  &= \pr{ 1 - \frac{\theta_{K-(k+2)}^4}{\theta_{K-(k+1)}^4} } (z_{k+1} - x_k^{+}).
\end{align}
Dividing left hand side with $h$ we observe,
\begin{align*}  
    \frac{x_{k+1}-x_{k}^+}{h}
        = \frac{x_{k+1}-x_{k} + h^2 \nabla f(x_{k})}{h}
        = \frac{x_{k+1}-x_{k}}{h} + O(h)   
        = \dot{X}(t) + O(h).
\end{align*}
Then from the fact $\theta_{K-k} = \frac{K-k}{2} + o(K)$, we observe
\begin{align*}
    \frac{1}{h} \pr{ 1 - \frac{\theta_{K-(k+2)}^4}{\theta_{K-(k+1)}^4} }    
        &= \frac{1}{h} \pr{ 1 - \frac{(K-k-2+o(K))^4}{(K-k-1+o(K))^4} }  \\
        &= \frac{1}{h} \pr{ \frac{(2 (K-k) + o(K)) ( 2 (K-k)^2-6 (K-k) + o(K)K)}{(K-k)^4 + o(K)K^3} } \\
        &= \frac{(2 (Kh-kh) + o(K)h) ( 2 (Kh-kh)^2 - 6 (Kh-kh)h + o(K)Kh^2)}{(Kh-kh)^4 + o(K)K^3h^4}    \\
        &= \frac{(2 (T-t) + o(K)h) ( 2 (T-t)^2 - 6 (T-t)h + o(K)Th)}{(T-t)^4 + o(K)T^3h}    
        = \frac{4}{T-t} + o(K)h.
\end{align*}
Dividing \eqref{eq : rewritten OGM-G x update} by $h$, applying above observations, 
corresponding $z_{k+1}$ with $Z(t+h)=Z(t)+O(h)$ we have
\begin{align*}
    \dot{X}(t) + O(h)
        = \frac{x_{k+1}-x_{k}^+}{h}
        = \frac{1}{h} \pr{ 1 - \frac{\theta_{K-(k+2)}^4}{\theta_{K-(k+1)}^4} } (z_{k+1} - x_{k}^{+})
        = \frac{4}{T-t}\pr{ Z(t) - X(t) } + O(h) + o(K)h.
\end{align*}
Organizing with respect to $Z$, 
we have
\begin{align*}
    Z(t) = \frac{T-t}{4} \dot{X}(t) + X(t) + O(h) + o(K)h.
\end{align*}
Now to conclude the desired result, we observe the followings.
First, observe the terms with gradient are $O(h)$.
For example, $\frac{1}{2L}\norm{\nabla f(x_K)}^2=\frac{h^2}{2}\norm{\nabla f(x_K)}^2 = O(h)$.
With this observation, we see $x_{k-1}^{+} = x_{k-1} - \frac{1}{L}\nabla f(x_{k-1})$ can be replaced with $x_{k-1}$.
Second, observe $h\theta_{K-k} = \frac{T-t}{2} + o(K)h$.
Third, we correspond $x_{k-1}$ with $X(t-h) = X(t) + O(h)$.

Plugging these to \eqref{eq:Discrete OGM-G Lyapunov}, and dividing by $2h^2$, we get
\begin{align*}
    &\frac{U_k}{2h^2}   
        = \frac{1}{2(h\theta_{K-k})^2} \pr{ f(x_k) - f(x_K) + O(h) } 
            + \frac{1}{2(h\theta_{K-k})^4} \inner{z_k-x_k+O(h)}{z_k-x_K+O(h)}    \\
        &= \frac{2}{(T-t + o(K)h)^2} \pr{ f(X(t)) - f(X(T)) } 
            + \frac{8}{(T-t + o(K)h)^4} \inner{Z(t)-X(t)}{Z(t)-X(T)}    +  O(h)        \\ 
        &= \frac{2}{(T-t)^2} \pr{ f(X(t)) - f(X(T)) } 
            + \frac{1}{2(T-t)^4} \inner{(T-t)\dot{X}(t)}{(T-t)\dot{X}(t) + 4(X(t)-X(T)) }  + O(h) + o(K)h         \\
        &= \frac{2}{(T-t)^2} \pr{ f(X(t)) - f(X(T)) } 
            + \frac{1}{2(T-t)^4} \pr{ \norm{(T-t)\dot{X}(t)+2(X(t)-X(T))}^2 - 4\norm{X(t)-X(T)}^2 } \\ &\quad
            + O(h) + o(K)h.
\end{align*}
Ignoring $O(h)$ and $o(K)h$, we see $\frac{U_k}{2h^2}$ corresponds to the Lyapunov function defined in Theorem \ref{thm:main_result_for_OGMG}.

\subsection{Details for Theorem \ref{thm:generalized_result_for_OGMG}} \label{appendix : Detail for generalized_result_for_OGMG}
Recall by plugging $\alpha=-2$, $X_c=X(T)$, $t_0=0$ to \eqref{eq : Energy for OGM-G for general r, alpha}, 
we obtained the conservation law presented in \ref{section : OGM-G r<-3}.
\begin{align*}   
    E   &\equiv \frac{2}{T^2} ( f(X_0) - f(X(T)) ) + \frac{r+3}{T^4} \norm{X_0-X(T)}^2 \\
        &= \frac{2}{(T-t)^{2}} \pr{ f(X) - f(X(T)) }  + \frac{1}{2(T-t)^{4}} \norm{ (T-t) \dot{X} + 2(X-X(T))}^2 
        + \frac{r+1}{(T-t)^{4}} \norm{X-X(T)}^2 \\&\quad
        + \int_{0}^{t} \frac{(-(r+3))}{(T-s)^{5}} \norm{(T-s) \dot{X} + 2(X-X(T))}^2 ds  \\&\quad
        + \int_{0}^{t} \frac{4}{(T-s)^{3}} \pr{  f(X(T)) - f (X) - \inner{ \nabla f (X) } { X(T) - X } }  ds .
\end{align*}
By collecting first three terms, define the Lyapunov function as
\begin{align*}
    \Phi(t) 
        &= \frac{2}{(T-t)^{2}} \pr{ f(X) - f(X(T)) } 
            + \frac{1}{2(T-t)^{4}} \norm{ (T-t) \dot{X} + 2(X-X(T))}^2
            + \frac{r+1}{(T-t)^{4}} \norm{X-X(T)}^2 .
\end{align*}
From conservation law we know $\dot{E}=0$, so we have
\begin{align*}
    \dot{\Phi}(t) = \frac{r+3}{(T-t)^{5}} \norm{(T-t) \dot{X} + 2(X-X(T))}^2 
                    - \frac{4}{(T-t)^{3}} \pr{  f (X(T)) - f(X) - \inner{ \nabla f (X) } { X - X(T) } } 
                    \le 0 .
\end{align*}
Note the first term is nonpositive since $r\le-3$.
Especially $\Phi(0)\ge \lim_{t\to T^{-}} \Phi(t)$.

Now we calculate $\lim_{t\to T^{-}} \Phi(t)$.
From \ref{appendix : Property of OGM-G ODE} we know $\lim_{t\to T^{-}} \frac{\dot{X}(t)}{t-T} = -\frac{2}{r+1} \nabla f(X(T))$.
By applying L'H\^{o}pital's rule we have
\begin{align*}
    \lim_{t\to T^{-}} \frac{f(X(t))-f(X(T))}{(T-t)^2}
        &= \lim_{t\to T^{-}} \frac{\inner{\nabla f(X(t))}{\dot{X}(t)}}{-2(T-t)}
        =  \inner{\nabla f(X(T))}{\lim_{t\to T^{-}} \frac{\dot{X}(t)}{2(t-T)} }        
        = -\frac{1}{r+1} \norm{\nabla f(X(T))}^2  \\
    \lim_{t\to T^{-}} \frac{X(t)-X(T)}{(T-t)^2}
        &= \lim_{t\to T^{-}} \frac{\dot{X}(t)}{-2(T-t)}
        = \frac{1}{2} \lim_{t\to T^{-}} \frac{\dot{X}(t)}{t-T}
        = -\frac{1}{r+1} \nabla f(X(T))  .
\end{align*}
Therefore we get
\begin{align*}
    \lim_{t\to T^{-}} \Phi(t)
        &= \lim_{t\to T^{-}} \pr{ \frac{2 \big( f(X)-f(X(T)) \big) }{(T-t)^2}
            + \frac{1}{2} \norm{ -\frac{\dot{X}}{t-T} + 2\frac{X-X(T)}{(T-t)^2} }^2 
            + (r+1) \norm{ \frac{X-X(T)}{(T-t)^2} }^2 } \\
        &= -\frac{2}{r+1} \norm{\nabla f(X(T))}^2 
            + \frac{1}{2} \norm{ \frac{2}{r+1} \nabla f(X(T)) -  \frac{2}{r+1} \nabla f(X(T)) }^2  
            + \frac{1}{r+1} \norm{\nabla f(X(T))}^2   \\
        &= \frac{1}{-(r+1)} \norm{\nabla f(X(T)) }^2 .
\end{align*}
Finally applying above calculation we have
\begin{align*}
    \frac{1}{-(r+1)} \norm{\nabla f(X(T)) }^2 
        = \lim_{t\to T^{-}} \Phi(t)    
        &\le \Phi(0)    
        = \frac{2}{T^2} ( f(X_0) - f(X(T)) ) + \frac{r+3}{T^4} \norm{X_0 - X(T)}^2  \\
        &\le \frac{2}{T^2} \pr{ f(X_0) - f(X(T)) }.
\end{align*}
This proves Theorem~\ref{thm:generalized_result_for_OGMG}.

\section{Proof of Theorem 5.1}       \label{appendix:proof_for_thm_5.1}
Recall, with $\theta_k = \frac{k}{2}$ the discretized method was
\begin{align}       \tag{\ref{eq:AGM ODE discretization}}
    x_{k}^+ &= x_k - \frac{s}{2} \nabla f(x_k)  \nonumber \\
    z_{k+1} &= z_{k} - s\theta_k \nabla f(x_k)    \nonumber \\ 
    x_{k+1} &= \frac{\theta_{k}^2}{\theta_{k+1}^2} x_{k}^{+} 
                + \pr{1-\frac{\theta_{k}^2}{\theta_{k+1}^2}} z_{k+1},       \nonumber 
\end{align}
and with $c_k=\frac{\theta_{k+1}}{\theta_{k+1}^2-\theta_k^2}$ the Lyapunov function was
\begin{align*}
   \Phi_k &= 2 c_k \theta_k^2 \pr{ f(x_k) - f_\star - \frac{s}{4} \norm{\nabla f(x_k) }^2 }
        +  \frac{1}{s} \norm{ z_{k+1} - X_\star }^2
\end{align*}
for $k=0,1,\dots$.
We first prove $\Phi_{k+1} \le \Phi_{k}$, then we will get the desired result from $\Phi_k \le \Phi_0$.

\textbf{(i) $\Phi_{k+1} \le \Phi_{k}$}    \\ 
For convenience, name
    $$ A_k = c_k \theta_k^2 = \frac{ \theta_{k+1} }{ \theta_{k+1}^2 - \theta_{k}^2 } \theta_k^2. $$
Observe since $c_k = \frac{2(k+1)}{(k+1)^2-k^2} = \frac{2(k+1)}{2k+1} \ge 1$, we have $A_k \ge \theta_k^2$.
From this we have
\begin{align*}
    \frac{1}{s} \norm{ z_{k+1} - X_\star }^2 - \frac{1}{s} \norm{ z_{k+2} - X_\star }^2    
    &= 2\theta_{k+1} \inner{\nabla f(x_{k+1})}{z_{k+1} - X_\star} 
        - s \theta_{k+1}^2 \norm{ \nabla f(x_{k+1}) }^2  \\
    &\ge 2\theta_{k+1} \inner{\nabla f(x_{k+1})}{z_{k+1} - X_\star}
        - s A_{k+1} \norm{ \nabla f(x_{k+1}) }^2 .
\end{align*}

Applying this fact we have
\begingroup
\allowdisplaybreaks
\begin{align*}
    \Phi_k-\Phi_{k+1}     
        &=  2 A_k \pr{ f(x_k) - f_\star - \frac{s}{4} \norm{\nabla f(x_k) }^2 } 
            - 2  A_{k+1} \pr{ f(x_{k+1}) - f_\star - \frac{s}{4} \norm{\nabla f(x_{k+1}) }^2 }        \\
            &\quad
            + \frac{1}{s} \norm{ z_{k+1} - X_\star }^2
            - \frac{1}{s} \norm{ z_{k+2} - X_\star }^2           \\
        &\ge  2 A_k \pr{ f(x_k) - f_\star - \frac{s}{4} \norm{\nabla f(x_k) }^2 } 
            - 2  A_{k+1} \pr{ f(x_{k+1}) - f_\star - \frac{s}{4} \norm{\nabla f(x_{k+1}) }^2 }        \\
            &\quad
            + 2\theta_{k+1} \inner{\nabla f(x_{k+1})}{z_{k+1} - X_\star}
            - s A_{k+1} \norm{ \nabla f(x_{k+1}) }^2          \\
        &=  2 A_k \pr{ f(x_k) - f_\star - \frac{s}{4} \norm{\nabla f(x_k) }^2 } 
            - 2  A_{k+1} \pr{ f(x_{k+1}) - f_\star + \frac{s}{4} \norm{\nabla f(x_{k+1}) }^2 }        \\
            &\quad
            + 2\theta_{k+1} \inner{\nabla f(x_{k+1})}{z_{k+1} - X_\star}    \\
        &=  2 A_k \pr{ f(x_k) - f_\star - \frac{s}{4} \norm{\nabla f(x_k) }^2 } 
            - 2  A_{k} \pr{ f(x_{k+1}) - f_\star + \frac{s}{4} \norm{\nabla f(x_{k+1}) }^2 }        \\
            &\quad
            + 2  \underbrace{\pr{ A_{k} - A_{k+1} + \theta_{k+1}} }_{=\frac{(k+1)^2}{8 k^2+16 k+6} \ge 0}  
                \pr{ f(x_{k+1}) - f_\star + \frac{s}{4} \norm{\nabla f(x_{k+1}) }^2 } \\
            &\quad
            - 2\theta_{k+1} \pr{ f(x_{k+1}) - f_\star + \frac{s}{4} \norm{\nabla f(x_{k+1}) }^2 }
            + 2\theta_{k+1} \inner{\nabla f(x_{k+1})}{z_{k+1} - X_\star}   \\
        &\ge  2 A_k \pr{ f(x_k) - f_\star - \frac{s}{4} \norm{\nabla f(x_k) }^2 } 
            - 2  A_{k} \pr{ f(x_{k+1}) - f_\star + \frac{s}{4} \norm{\nabla f(x_{k+1}) }^2 }        \\
            &\quad
            - 2\theta_{k+1} \pr{ f(x_{k+1}) - f_\star + \frac{s}{4} \norm{\nabla f(x_{k+1}) }^2 }
            + 2\theta_{k+1} \inner{\nabla f(x_{k+1})}{z_{k+1} - X_\star}   \\
        &= 2 A_k 
            \pr{ f(x_k) - f(x_{k+1}) - \frac{s}{4} \norm{\nabla f(x_k) }^2 - \frac{s}{4} \norm{\nabla f(x_{k+1}) }^2}     \\  
            &\quad
            + 2\theta_{k+1} \pr{ f_\star - f(x_{k+1}) - \inner{\nabla f(x_{k+1})}{X_\star-x_{k+1}} - \frac{s}{4} \norm{\nabla f(x_{k+1}) }^2  }        \\  &\quad
            + 2\theta_{k+1} \inner{\nabla f(x_{k+1})}{z_{k+1} - x_{k+1}}   \\
        &\stackrel{(a)}{\ge} 2 A_k 
            \pr{ f(x_k) - f(x_{k+1}) - \frac{s}{4} \norm{\nabla f(x_k) }^2 - \frac{s}{4} \norm{\nabla f(x_{k+1}) }^2}     \\  
            &\quad
            + 2\theta_{k+1} \inner{\nabla f(x_{k+1})}{z_{k+1} - x_{k+1}}   \\
        &= 2 A_k 
            \pr{ f(x_k) - f(x_{k+1}) - \frac{s}{4} \norm{\nabla f(x_k) }^2 - \frac{s}{4} \norm{\nabla f(x_{k+1}) }^2}     \\  
            &\quad
            + 2\theta_{k+1} \inner{\nabla f(x_{k+1})}{ \frac{ \theta^2_{k} }{ \theta^2_{k+1} - \theta^2_{k} } \pr{ x_{k+1} - x_{k}^+ }}   \\
        &= 2 A_k 
            \pr{ f(x_k) - f(x_{k+1}) - \frac{s}{4} \norm{\nabla f(x_k) }^2 - \frac{s}{4} \norm{\nabla f(x_{k+1}) }^2}     \\  &\quad
            + 2 A_k \inner{\nabla f(x_{k+1})}{ x_{k+1} - x_k + \frac{s}{2} \nabla f (x_k) }   \\
        &= 2 A_k 
            \pr{ f(x_k) - f(x_{k+1}) + \inner{\nabla f(x_{k+1})}{ x_{k+1} - x_k } - \frac{s}{4} \norm{ \nabla f(x_k) - \nabla f(x_{k+1}) }^2 }  
        \stackrel{(b)}{\ge} 0 .
\end{align*}
\endgroup
The inequalities $(a)$ and $(b)$ come from the fact $s\in\left( 0, \frac{2}{L} \right]$ and $L$-smoothness of $f$.

\textbf{(ii) From $\Phi_{k} \le \Phi_{0}$, we have $f(x_k^{+}) - f_\star \le  \frac{k+\frac{1}{2}}{k+1} \frac{2\norm{ X_0 - X_\star }^2}{s k^2} $ }    \\ 
From $\theta_0=0$ we have $A_0=0$, and so $z_1=z_0+s\theta_0\nabla f(X_0) = z_0 = X_0$.
Therefore
\[
    \Phi_0 = 2A_0 + \frac{1}{s} \norm{z_1-X_\star}^2 
         = \frac{1}{s} \norm{X_0-X_\star}^2
\]
Now since $f$ is $L$-smooth, for $s\in\left(0,\frac{2}{L}\right]$, we have 
\[
    f(x_k^+)    
        \le  f(x_k) - \frac{1}{2L} \norm{\nabla f(x_k) }^2
        \le  f(x_k) - \frac{s}{4} \norm{\nabla f(x_k) }^2,
\]
and so
\begin{align*}
    2 A_k \pr{ f(x_k^+)-f_\star }  
        &\le 2 A_k \pr{ f(x_k) - f_\star - \frac{s}{4} \norm{\nabla f(x_k) }^2 }   
        \le \Phi_k   
        \le \Phi_0 = \frac{1}{s} \norm{ X_0 - X_\star }^2.
\end{align*}
Therefore, we conclude
\begin{align*}
    f(x_k^+)-f_\star 
        &\le \frac{\norm{ X_0 - X_\star }^2}{ 2s A_k }  
        = \pr{ \frac{ \theta_{k+1} }{ \theta_{k+1}^2 - \theta_{k}^2 } \theta_k^2  }^{-1} \frac{\norm{ X_0 - X_\star }^2}{ 2s }  \\
        &= \pr{ \frac{ 2k + 1 }{ 2(k + 1) } \times \frac{4}{k^2}  }\frac{\norm{ X_0 - X_\star }^2}{ 2s }    \\
        &= \frac{ k+\frac{1}{2} }{ k+1 } \frac{2\norm{ X_0 - X_\star }^2}{ sk^2 }    .
\end{align*}
Since $\frac{ k+\frac{1}{2} }{ k+1 }\le1$, this implies $f(x_k^+)-f_\star \le \frac{2\norm{ X_0 - X_\star }^2}{ sk^2 }$ as well.
This proves Theorem~\ref{thm: convergence rate for r=3}.
\qed

\section{Time-dependent Hamiltonian}
\label{s:hamiltonian}
For the sake of completeness, we show how the dynamics is described through a Hamiltonian perspective.
With the Hamiltonian
\begin{align*}
H(W,P,t)&=\langle P, \dot{W}\rangle-L(W,P,t)\\
&=
\frac{t}{2}\|P\|^2+t^3(f(X(W,t))-f_\star),
\end{align*}
the dynamics of the Euler--Lagrange equation can be equivalently specified with
\begin{align*}
    \dot{P} &= -\nabla_W H(W,P,t) 
            = -t\nabla f(X(W,t))   \\
    \dot{W} &= \nabla_PH(W,P,t)  
            = tP.
\end{align*}
However, our setup differs from the classical setup in that the Lagrangian and the Hamiltonian explicitly depend on time.
One consequence of this difference is that the Hamiltonian is not conserved:
\begin{align*}
    \frac{d}{dt} H(W,P,t) 
    &= \inner{\dot{W} }{\nabla_W H(W,P,t)}
       + \inner{ \dot{P} }{ \nabla_P H(W,P,t) }
       + \pd{t} H (W,P,t)     \\
    &= \inner{ \nabla_P H(W,P,t) }{\nabla_W H(W,P,t) }
       + \inner{ -\nabla_W H(W,P,t) }{ \nabla_P H(W,P,t) }
       + \pd{t} H (W,P,t)     \\&
    = \pd{t} H (W,P,t) \ne 0.
\end{align*}
Since $H$ is not conserved, the classical theory of symplectic integrators is not immediately applicable.

\end{document}